\documentclass[letterpaper, 10 pt, journal, twoside]{ieeetran}

\IEEEoverridecommandlockouts                              

\usepackage{graphics} 
\usepackage{epsfig} 
\usepackage{epstopdf}
\usepackage{amsmath} 
\usepackage{amssymb}  
\usepackage{amsfonts}
\usepackage{subfigure}
\usepackage{nicefrac}
\usepackage{cleveref}
\usepackage{mathtools}
\usepackage{framed}
\usepackage{tabularx} 

\DeclarePairedDelimiter{\ceil}{\lceil}{\rceil}
\DeclarePairedDelimiter{\floor}{\lfloor}{\rfloor}

\usepackage[T3,T1]{fontenc}

\DeclareSymbolFont{tipa}{T3}{cmr}{m}{n}
\DeclareMathAccent{\invbreve}{\mathalpha}{tipa}{16}

\usepackage{xspace}
\usepackage{soul}
\usepackage{mathtools}
\usepackage{makecell}
\usepackage{todonotes}
\usepackage{nicefrac}

\newcommand{\ra}{\rightarrow}

\newcommand{\ie}{\unskip, i.\,e.,\xspace}
\newcommand{\eg}{\unskip, e.\,g.,\xspace}

\newcommand{\sut}{\text{s.\,t.\,}}

\newcommand{\nrm}[1]{\left\lVert#1\right\rVert}

\newcommand{\scal}[1]{\left\langle#1\right\rangle}
\newcommand{\tr}[1]{\ensuremath{\text{tr}}\left(#1\right)}
\newcommand{\E}[1]{\mathbb E\left[#1\right]}

\newcommand{\N}{\ensuremath{\mathbb{N}}}

\newcommand{\R}{\ensuremath{\mathbb{R}}}

\newcommand{\X}{\ensuremath{\mathbb{X}}}

\newcommand{\U}{\ensuremath{\mathbb{U}}}

\newcommand*\diff{\mathop{}\!\mathrm{d}}
\newcommand{\eps}{\ensuremath{\varepsilon}}

\newcommand{\spc}{\ensuremath{\,\,}}

\newcommand{\dom}{\ensuremath{\text{dom}}}

\newcommand{\ball}{\ensuremath{\mathcal B}}

\newcommand{\co}{\ensuremath{\overline{\text{co}}}}

\DeclareMathOperator*{\argmin}{arg\,min}

\definecolor{dgreen}{rgb}{0.0, 0.5, 0.0}

\newtheorem{dfn}{Definition}

\newtheorem{lem}{Lemma}
\newtheorem{thm}{Theorem}
\newtheorem{crl}{Corollary}
\newtheorem{rem}{Remark}

\newcommand{\SH}{S\&H\xspace}

\newcommand{\D}{\ensuremath{\mathcal{D}}}
\newcommand{\A}{\ensuremath{\mathcal{A}}}
\newcommand{\K}{\ensuremath{\mathcal{K}}}

\newcommand{\lip}[2]{\ensuremath{\text{Lip}_{#2}\left(#1\right)}}

\newcommand{\Var}[1]{\mathbb V\left[#1\right]}

\newcommand{\ao}[1]{\alpha_1(#1)}
\newcommand{\at}[1]{\alpha_2(#1)}
\newcommand{\aao}[1]{{^\lambda\alpha_1(#1)}}
\newcommand{\aat}[1]{{^\lambda\alpha_2(#1)}}
\newcommand{\aaoi}[1]{{^\lambda\alpha^{-1}_1(#1)}}

\newcommand{\aaoc}[1]{{^\lambda\breve{\alpha}_1(#1)}}
\newcommand{\aatc}[1]{{^\lambda\invbreve{\alpha}_2(#1)}}
\newcommand{\aaoci}[1]{{^\lambda\breve{\alpha}^{-1}_1(#1)}}

\newcommand{\adci}[1]{{\breve{\alpha}^{-1}_3(#1)}}
\newcommand{\adi}[1]{{\alpha^{-1}_3(#1)}}
\newcommand{\adc}[1]{{\breve{\alpha}_3(#1)}}
\newcommand{\ad}[1]{{\alpha_3(#1)}}

\newcommand{\nm}[1]{\lVert #1\rVert_2}
\newcommand{\intstep}[1]{\int_{k\delta}^{(k + 1)\delta} #1 \diff t}
\newcommand{\intstepg}[1]{\int_{k\delta}^{k\delta + t'} #1 \diff t}

\newcommand{\rhs}{f(X_t, U_t) + \sigma(X_t, U_t)Z_t}
\newcommand{\fb}{\bar{f}}
\newcommand{\sgb}{\bar{\sigma}}
\newcommand{\Zb}{\bar{Z}}

\newcommand{\rtd}{\tilde{r}}
\newcommand{\rht}{\hat{r}}
\newcommand{\rb}{\bar{r}}
\newcommand{\rst}{r^*}
\newcommand{\Rst}{R^*}

\renewcommand{\eps}{\varepsilon}
\renewcommand{\P}[1]{\mathbb{P}\left[#1\right]}

\newcommand{\La}[1]{L_\lambda\left(#1\right)}
\newcommand{\evb}{\bar{\mu}}
\newcommand{\sdb}{\tilde{\sigma}}
\newcommand{\dprod}[2]{\langle #1, #2 \rangle}

\newcommand{\Lip}{\text{Lip}}

\newcommand{\infconv}[1]{\inf\limits_{y \in \mathbb{X}}\left(#1(y) + \frac{\nrm{y - x}^2}{2\lambda^2}\right)}
\newcommand{\al}[1]{\begin{aligned} #1 \end{aligned}}

\newcommand{\tsqrt}{\lambda\sqrt{2\at{R}}}

\title{On stochastic stabilization via non-smooth control Lyapunov functions}

\author{Pavel Osinenko, Grigory Yaremenko, Georgiy Malaniya
\thanks{The authors are with Skolkovo Institute of Science and Technology, emails: \texttt{\{p.osinenko,grigory.yaremenko,g.malaniya\} @skoltech.ru}.
The experimental case study of this work was supported by Dmitry Dobriborsci and Ksenia Makarova.
}}
\begin{document}

\maketitle
\thispagestyle{empty}
\pagestyle{empty}

\begin{abstract}
Control Lyapunov function is a central tool in stabilization.
It generalizes an abstract energy function -- a Lyapunov function -- to the case of controlled systems.
It is a known fact that most control Lyapunov functions are non-smooth -- so is the case in non-holonomic systems, like wheeled robots and cars.
Frameworks for stabilization using non-smooth control Lyapunov functions exist, like Dini aiming and steepest descent.
This work generalizes the related results to the stochastic case.
As the groundwork, sampled control scheme is chosen in which control actions are computed at discrete moments in time using discrete measurements of the system state.
In such a setup, special attention should be paid to the sample-to-sample behavior of the control Lyapunov function.
A particular challenge here is a random noise acting on the system.
The central result of this work is a theorem that states, roughly, that if there is a, generally non-smooth, control Lyapunov function, the given stochastic dynamical system can be practically stabilized in the sample-and-hold mode meaning that the control actions are held constant within sampling time steps.
A particular control method chosen is based on Moreau-Yosida regularization, in other words, inf-convolution of the control Lyapunov function, but the overall framework is extendable to further control schemes.
It is assumed that the system noise be bounded almost surely, although the case of unbounded noise is briefly addressed.
\end{abstract}


\section{Introduction}
\label{sec_intro}

Stochastic stability theory can be traced back to the works of Khasminskii \cite{khasminskii201-stochastic}, Kushner \cite{Kushner1965stochastic-stability}, Mao \cite{mao1991-stability} and Deng et al. \cite{Deng2001-stochastic-stab-noise} who extended the classical results and translated them into the language of $\K$-functions of Khalil \cite{Khalil1996-nonlin-sys}, convenient for control engineers. 
Stochastic stability analyses were applied to various types of systems including discrete systems \cite{McAllister2003-stochastic}, cascaded systems \cite{Liu2008-stochastic}, delayed systems \cite{Liu2011-stochastic}, systems with input saturation \cite{Li2017-Stochastic}, systems with state-dependent switching \cite{Wu2013-stochastic}, nonlinear stochastic dynamic systems with singular perturbations \cite{Cheng2001-stochastic-stability}, hybrid systems \cite{Teel2014-stochastic-stability}, hybrid retarded systems \cite{Huang2008-hybrid-stochastic-retarded}, linear systems with randomly jumping parameters \cite{Feng1992-stochastic-stability} etc.
Practical stochastic stability was addressed in \eg \cite{caraballo2015-practical-stochastic,Qin2020-stochastic,Do2020-stochastic}.
When it comes to stochastic stabilization, most of the existing works assume continuous application of the control, whereas sampled control schemes were considered in rather specific contexts \eg based on approximate discrete-time models \cite{fu2016sampled,yu2018sampled} or in the event-triggered mode \cite{GAO2018}.
Stabilization of stochastic logical systems \cite{yang2022stability, li2020robustness} is yet another example of stochastic stabilization performed in discrete time, however the latter setting also implies a finite state-space, which is not suitable for certain applications \eg in robotics.

The importance of sampled control frameworks is dictated by two facts.
First, most modern controllers are implemented on digital media whence control actions are naturally computed at discrete moments in time, and based on discrete measurement of the state, rather than continuous.
Second, the need for sampled controls is motivated by the fact that the majority of dynamical systems are not even stabilizable by means of feedback laws that depend continuously on the state \cite{Brockett1983-stabilization,Sontag1990-stabilization-survey,Clarke1997-stabilization,Fontes2003discontinuous,Cortes2008-discont-dyn-sys,Clarke2009-slid-mode-stab,Braun2017-SH-stabilization-Dini-aim}.
There is an intimate connection between discontinuous feedback controls and non-smooth control Lyapunov functions \cite{Clarke2011-discont-stabilization}.
In particular, non-existence of a smooth control Lyapunov function leads to non-existence of a continuous stabilizing feedback law.
It should be noted here that not only do control Lyapunov functions appear mostly non-smooth, numerical routines for their calculation commonly produce non-smooth functions -- for particular methods, refer \eg to \cite{Baier2012-linear,Baier2014-num-CLF,grammatico2013control,Giesl2015-review}.
A major problem that arises under discontinuous feedback laws is defining the system trajectory.
But even resorting to such generalized notions of solutions as of Filippov does not remedy the situation.
It was the prominent work by Clarke, Ledyaev, Sontag and Subbotin \cite{Clarke1997-stabilization} which derived a general framework for stabilization by means of non-smooth control Lyapunov functions by means of sampled controls.
Stabilization here was meant as practical stabilization \ie the system state could be stabilized into any desirable vicinity of the equilibrium, provided that the sampling time be sufficiently small.
In turn, the control law computation was based on the Moreau-Yosida regularization, in other words, inf-convolution of the control Lyapunov function involved.
Further methods are applicable -- for surveys, refer \eg to \cite{Braun2017-SH-stabilization-Dini-aim,Osinenko2019nonsmoothstabsurvey}.

A brief remark should be made on the case of smooth control Lyapunov functions.
Here, one can exploit universal formulas for stabilizing controls \cite{Sontag1989-formula,Lin1991-stabilization-bounded-controls}.
Their generalizations to the stochastic case exist.
For instance, Florchinger \cite{Florchinger1994-Lyapunov-like-stochastic-stability} suggested constructions of stochastically stabilizing controls via Lie theory methods.
Smooth control Lyapunov functions were also assumed in the stochastic stabilization work of Deng et al. \cite{Deng2001-stochastic-stab-noise} and Gao et al. \cite{GAO2018}, which addressed event-triggered digital implementation of a given feedback law that stabilizes the system in second moment.
In the present work, in contrast, non-smooth control Lyapunov functions are considered.

There is no extension of the Clarke's general framework for stabilization to the case of stochastic systems.
Whereas Deng et al. \cite{Deng2001-stochastic-stab-noise} considered the continuous control case, the work \cite{Osinenko2021stochSH} extended on these results in the case of sampled controls.
However, the control Lyapunov function involved was smooth, whereas the non-smooth case was tackled only briefly, specifically, under bounded noise.
The current work elaborates on the matter of stochastic stabilization via sampled controls computed from non-smooth control Lyapunov functions to a full technical extent.

The engineering background of the proposed framework for stochastic stabilization is motivated by the two facts.
First, as also mentioned above, most modern controllers are realized in digital devices that naturally lead to a sampled control scheme where control actions are computed at discrete, usually evenly distributed, moments in time upon receiving state measurements at discrete moments in time as well.
In the presented work, the control law is treated explicitly in the sample-and-hold mode, whereas the underlying system dynamics model is considered as time-continuous, to be precise, in the form of a stochastic differential equation.
Second, the formalism of stochastic differential equations allows accounting for random uncertainty that may be related to model imperfection, system, measurement and/or actuator noise etc.
Effects of noise magnitude and sampling time choice are demonstrated in an experimental case study with a mobile robot in Section \ref{sec_experiment}.

\textbf{Contribution.}
Informally, the contribution of the current work is summarized as follows.
In Theorem \ref{thm_bounded_noise} it is shown that if a given non-linear stochastic control system with bounded noise has a, generally non-smooth, control Lyapunov function, a sample-and-hold policy can be produced that stabilizes the system's equilibrium.
The theorem explicitly provides the aforementioned policy and estimates the quality of resulting stabilization for a given configuration of noise, sampling rate and computation error. Unlike previously mentioned results in stochastic stabilization, Theorem \ref{thm_bounded_noise} concerns the case of a potentially non-smooth control-Lyapunov function and at the same time accounts for imperfections peculiar to a digital controller, namely control update latency and computation error.

\begin{center}
\textbf{Table of notation} 
\hphantom{.}\\
\hphantom{.}\\
\begin{tabularx}{0.98\columnwidth}{ |p{1.5cm}|p{6.29cm}| }  
 	\hline
	$\X$ & State space $\R^n$, where $n \in \N$ \\ 
	$\U$ & Control set, which is a subset of $\R^m$, \\
	& where $m \in \N$ \\
	$\mathbb{R}_+$ & $[0, +\infty)$\\
	$\E{\cdot}$ & Expected value\\ 
	$\Var{\cdot}$ & Variance\\
	$\P{\cdot}$ & Probability measure\\
	$\mathbb{R}_+$ & $[0, +\infty)$\\
	$\nrm{\cdot}$ & Eucledian norm $\nrm{\cdot}_2$\\ 
	$\nrm{\cdot}_{\text{op}}$ & Operator norm\\ 
	$\K_\infty$ & Class of kappa-infinity functions \cite{Khalil1996-nonlin-sys}\\
	$\ball_R$ & Closed ball of radius $R$ centered \\
	& at the origin \\
	$\Lip_h(R)$ & Lipschitz constant of $h(\cdot)$ in $\ball_R$ \\  
	$a\land b$ & $\min(a, b)$ \\
	$\scal{\cdot, \cdot}$ & Scalar product \\
	$\text{tr}(\cdot)$ & Matrix trace \\
	$\mathcal{N}(x, y)$ & Normal distribution with mean $x$ and \\ & variance $y$ \\
	$\bar{\mathcal{N}}(x, y, z)$ & Truncated normal distribution with mean $x$, \\ & variance $y$ and maximal deviation $z$\\
	$\nabla\cdot$ & Gradient as a row-vector \\
 	\hline
\end{tabularx}
\begin{tabularx}{0.98\columnwidth}{ |p{1.5cm}|p{6.29cm}| } 
\hline  
  
	$\breve{h}$ & Lower convex envelope of function $h$ \\
	& $\breve{h} = \sup\{g \ | \ g\text{ is convex}, g(x) \leq f(x)\}$ \\
	$\invbreve{h}$ & Upper concave envelope of function $h$ \\
	& $\invbreve{h} = \inf\{g \ | \ g\text{ is concave}, g(x) \geq f(x)\}$ \\
	$\co(A)$ & Closed convex hull of set $A$\\
	$\mathcal{D}_vf(x)$ & Lower directional Dini derivative \\
	& \cite{Clarke1997-stabilization} of $f(\cdot)$ \\
	&  at $x$ in direction $v$ \\
	$\circ^{l}h(\cdot)$ & $l$-fold iteration of function $h(\cdot)$, \\
	&  like so $\underbrace{h(h(h(\dots)))}\limits_{l \text{ times}}$ \\
	$\mathbb I_A(\cdot)$ & Indicator function of set $A$ \\ 
	$a\cdot b$ & $a$ multiplied by $b$, where \\
	& $a$ and $b$ are scalars \\
	$\ceil{x}$ & $x$ rounded to the closest \\
	& greater or equal integer \\
	$\floor{x}$ & $x$ rounded to the closest \\
	& less or equal integer \\
	$a\text{ mod }b$ & remainder of $a$ when divided by $b$  \\

 \hline
\end{tabularx}
\end{center}

\textbf{Structure of the paper.}
All proofs and auxiliary lemmas are placed in the appendix.
The main body contains the technical preliminaries; the main theorem on stochastic stabilization by means of non-smooth control Lyapunov functions under almost surely bounded driving noise; the discussion of the case of unbounded driving noise; the demonstrative experimental study.


\section{Preliminaries}
\label{sec_prelim}

The aim of this work is to address stabilization of stochastic control systems of the class
\begin{equation}
\label{eqn:sys-core}
	\diff  X_t = f(X_t, U_t) \diff t + \sigma (X_t, U_t) Z_t \diff t,
\end{equation}
where $\{X_t\}_t, \{U_t\}_t$ are the state and, respectively, control stochastic processes, $\X$-, respectively, $\U$-valued; $f: \X \times \U \ra \R^n, \sigma: \X \times \U \ra \R^{n \times d}$; $\{Z_t\}_t$ is a $d$-dimensional random process that is measurable with respect to $t$.
The technical goal is to study stabilization of \eqref{eqn:sys-core} by Markov control policies in sampled mode. 
``Sampled mode" here means that the controller computes control actions based on measurements of state at discrete equidistant moments in time and keeps these actions constant between the said moments in time.
This is stated mathematically precise in Section \ref{sec_main_thm}.

Now, proceed to the necessary definitions for stability.
Consider a general stochastic system
\begin{equation}
	\label{eqn:sys-general}
	\begin{aligned}
		\diff  X_t & = f(X_t, t) \diff t + \sigma (X_t, t) \diff W_t, \\
	\end{aligned}
\end{equation}
where $\diff W_t$ is a placeholder for either $Z_t \diff t$ or $\diff B_t$ with $\{B_t\}_t$ being a standard Brownian motion.
The latter case will be used in the discussion on the unbounded noise case of Section \ref{sec_unbounded}.
Observe that if we were to apply a concrete policy $V_t$ to the control system \eqref{eqn:sys-core} by asserting $U_t := V_t$, we would get a system of the kind \eqref{eqn:sys-general} provided the driving noise is of the form $Z_t \diff t$.
Now, we are ready to state the key definitions for stability.

\begin{dfn}[Semi-asymptotic stability in probability]
	\label{dfn_conv}
	The origin of a stochastic system \eqref{eqn:sys-general} is said to be semi-asymptotically stable in probability in $R > 0$ until $r \geq 0$ if:
	\begin{equation}
		\label{eqn:conv}
		\begin{aligned}
			& x_0 \in \ball_R \implies \P{\limsup_{t \ra \infty} \nrm{X_t} \le r} = 1.
		\end{aligned}
	\end{equation}
\end{dfn}

This notion is similar to (almost sure) local asymptotic stability.
The only difference is that with asymptotic stability the attractor is a single point -- the origin, whereas with semi-asymptotic stability the attractor is a ball centered at the origin.
Thus $r$ is the radius of the attractor and $R$ is the radius of the basin of attraction.
If a stochastic system is said to be semi-asymptotically stable, this kind of attraction will occur almost surely.
The need to introduce this notion arises from how adding noise to an asymptotically stable deterministic system impacts the said system's stability.
Unless the magnitude of added noise is too large, the system will generally preserve its attractive properties within some basin of attraction ($\ball_R$), but the attraction will cease in some small neighborhood of the origin ($\ball_r$).

\begin{dfn}[Semi-asymptotic stability on average]
	\label{dfn_conv_mean}
	The origin of a stochastic system \eqref{eqn:sys-general} is said to be semi-asymptotically stable on average in $R > 0$ until $r \geq 0$:
	\begin{equation}
		\label{eqn:conv-mean}
		\begin{aligned}
			& x_0 \in \ball_R \implies \limsup_{t \ra \infty} \E{\nrm{X_t}} \le r.
		\end{aligned}
	\end{equation}
\end{dfn}

The latter definition is analogous to Definition \ref{dfn_conv}, but describes stability in mean.
Here the expected value of the state is guaranteed to get arbitrarily close to $\ball_r$ and stay there permanently, provided that the initial state was within $\ball_R$.

The following definitions are required to state and prove Theorem \ref{thm_bounded_noise}, which is the central result of this work.

\begin{dfn}
	\label{dfn:rgb}
	A non-negative locally Lipschitz-continuous function \\ $\rho \ : \ \mathbb{R}_+ \rightarrow \mathbb{R}_+$ is called a \textbf{radial growth bound} of system \eqref{eqn:sys-main} iff
	\begin{equation}
	\begin{aligned}
	& \rho(\nrm{x})\nrm{x} \geq x^Tf(x, u) + \nrm{\sigma^T(x, u)x}\bar{Z}, \ \forall u \in U. 
	\end{aligned}
	\end{equation}
	\end{dfn}
	
	\begin{dfn}
	\label{dfn:rf}
	A two-variable function $y(R;t)$ is called a \textbf{radial forecast} of \eqref{eqn:sys-core} iff there exists such a radial growth bound $\rho$ that for each $R \geq 0$, $y(R;\cdot)$ equals to the (local) solution of:
	\begin{equation}
	\label{eqn:radial-forecast}
	\begin{cases}
	\begin{aligned}
	 & \dot{y} = \rho(y), \\
	 & y(0) = R.
	\end{aligned}
	\end{cases}
	\end{equation}
\end{dfn}


\begin{rem}
	Note that for a strictly positive radial growth bound $\rho$, 
	$$
	y(R;t) = y(0; t + y^{-1}(0; R)),
	$$
	where the inversion $y^{-1}$ is with respect to the second argument.
\end{rem}

\begin{dfn}
	\label{dfn_CLF}
	A function $L$ is called a control Lyapunov function for \eqref{eqn:sys-core} if it satisfies:
	\begin{equation}
		\begin{aligned}
			& \forall x \in \X \spc \inf_{\nu \in \co(f(x, \U))} \D_\nu L(x) \le - \alpha_3(\nrm{x}),	\\
			&\alpha_1(\nrm{x}) \leq L(x) \leq \alpha_2(\nrm{x}), \spc \alpha_1, \alpha_2, \alpha_3 \in \K_\infty.
		\end{aligned}
	\end{equation}
\end{dfn}

\begin{dfn}
	Let there be a control Lyapunov function $L$ for \eqref{eqn:sys-core}
	and let $y_R$ be a radial forecast of \eqref{eqn:sys-core}, then a partial function $\tilde{r} \ : \ \mathbb{R}_+^5 \rightarrow \mathbb{R}$ is called an \textbf{attraction function} of \eqref{eqn:sys-core} if:
	 
	\begin{equation}
	\label{eqn:attraction-function}
		\al{
		&\tilde{r}(R, \Zb, \lambda, \delta, \eta) := \alpha_3^{-1}\Bigl( \eta + \\ 
		& + \frac{\sqrt{2\alpha_2(R)}}{\lambda}\Big(2\Lip_L\big(R + \\
		& + \lambda \sqrt{2\at{R}}\big)\Lip_f\big(R + \lambda \sqrt{2\at{R}}\big)\lambda^2 + \\ 
		& + \frac{\delta}{2}\Lip_f(y(R; \delta))(\fb(y(R;\delta)) + \sgb(y(R;\delta))\Zb) + \\ 
		& + \sgb(y(R;\delta))\Zb\Big) + \frac{\delta(\fb(y(R;\delta)) + \sgb(y(R;\delta))\Zb)^2}{2\lambda^2}\Bigr) + \\
		& + \lambda \sqrt{2\alpha_2(R)},
		}
	\end{equation}
	where for each $r \geq 0$ we have
	\begin{equation}
	\label{eqn:core-bounds}
		\begin{aligned}
		& \bar{f}(r) \geq \nrm{f(x, u)}, \forall x \in \mathcal{B}_r, \ u \in \mathbb{U},\\
		& \bar{\sigma}(r) \geq \nrm{\sigma(x, u)}_{op}, \forall x \in \mathcal{B}_r, \ u \in \mathbb{U}, \\
		& \nrm{f(x, u) - f(y, u)} \leq \text{Lip}_f(r)\nrm{x - y}, \ \forall x, y \in \mathcal{B}_r, u \in \mathbb{U},\\
		& \nrm{L(x) - L(y)} \leq \text{Lip}_L(r)\nrm{x - y}, \ \forall x, y \in \mathcal{B}_r, \\
		& \text{Lip}_L(\cdot), \text{Lip}_f(\cdot), \bar{f}(\cdot), \bar{\sigma}(\cdot) \text{ \, -- are non-decreasing and} \\ & \text{upper-semicontinuous functions.} 
		\end{aligned}
	\end{equation} 
\end{dfn}
For brevity, we denote $\rht_{\evb, \Zb, \delta, \eta}(r) := \rht(r, \Zb, \lambda, \delta, \eta)$.

We define $\La{\cdot}$ as the infimal convolution with a parameter $\lambda >0$ of $L(\cdot)$ as follows:
\begin{equation}
\al{
& \La{x} := \infconv{L}, \\
& \aao{\nrm{x}} \leq \La{x} \leq \aat{\nrm{x}}, \spc ^\lambda \alpha_1, ^\lambda \alpha_2 \in \K_\infty.
}
\end{equation}
Here, we used the upper left index $\lambda$ for notation purposes to stress the relation to the infimal convolution with a parameter $\lambda$.
One way to obtain $\aao{\cdot}$ and $\aat{\cdot}$ is to compute the infimal convolutions of $\ao{\cdot}$ and $\at{\cdot}$.
\begin{dfn}
A function $\rht \ : \ \mathbb{R}_+^6 \rightarrow \mathbb{R}$ is called a \textbf{mean attraction function} of \eqref{eqn:sys-core} if:
\begin{equation}
	\al{
	&\rht(r, \evb, \sdb, \lambda, \delta, \eta) := \breve{\alpha}_3^{-1}\Bigl(\eta + \\
	& + \sqrt{2\at{r}}\Big(2 \lambda \Lip_L(r + \lambda \sqrt{2\at{r}})\Lip_f(r + \\
	& + \lambda \sqrt{2\at{r}}) + \frac{\delta}{2a}\Lip_f(r)(\fb(r) +\\
	& + \sgb(r)\evb) + \frac{\sgb(r)\evb}{\lambda}\Big) + \\
	& + \frac{2\delta}{\lambda^2}((\fb(r) + \sgb(r)\Zb)^2 + \sdb^2\sgb(r)^2)\Bigr) + \\
	& + \lambda \sqrt{2\alpha_2(r)},
	}
\end{equation}

\end{dfn}
where $\adc{\cdot}$ is the lower convex envelope of $\alpha_3$ on $[0, r]$.
For brevity, let us denote $\rht_{\evb, \sdb, \lambda, \delta, \eta}(r) := \rht(r, \evb, \sdb, \lambda, \delta, \eta)$.

Naturally, an attraction function may only exist if $L(x)$ and $f(x, u)$ are locally Lipschitz continuous w.r.t. $x$. 
Also note that if $\text{Lip}_L(\cdot), \text{Lip}_f(\cdot), \bar{f}(\cdot), \bar{\sigma}(\cdot)$ fail to be non-decreasing, it would be easy to construct alternative bounds that are not only non-decreasing, but are in fact sharper than the original ones (i.e. by evaluating $\hat{f}(r) = \inf\limits_{r \leq r'}\bar{f}(r')$).

In Lemma \ref{lem:af-maf} it is proven that both the attraction and mean attraction functions exist whenever bounds in \eqref{eqn:core-bounds} exist.

\begin{rem}
	The domain of definition of $\tilde{r}(\cdot, \Zb, \lambda, \delta, \eta)$ depends on $\delta$, however for each $R$, $\tilde{r}(R, \Zb, \lambda, \delta, \eta)$ exists, provided that $\delta$ is sufficiently small. Unlike $\rtd$, the domain of definition of $\rht$ is not restricted.  
\end{rem}

Additionally we assume that some bounds for moments of $\nrm{Z_t}$ are known, in particular, there exist numbers $\bar \mu, \tilde \sigma$ such that:
\begin{equation}
	\al{
	&\forall t \spc \evb \geq \E{\nrm{Z_t}}, \\
	&\forall t \spc \sdb^2 \geq \Var{\nrm{Z_t}}.
	}
\end{equation}

\section{Main theorem}
\label{sec_main_thm}

Sample-and-hold mode introduces a latency between a change in the systems state and the controller's response to that change.
Considering such a setting yields a practical advantage: unlike feedback of the kind $U_t = \mu(X_t)$, sample-and-hold mode accurately describes policies that digital controllers can implement.
Furthermore, it provides other benefits.
A system produced by asserting $U_t := \mu(X_t)$ may have no solutions \eg
\begin{equation}
	\begin{cases}
	& \diff X_t = U_t \diff t, X_0 = 1,\\
	& U_t := 1 - 2 \mathbb I_{\mathbb{R}_+}(X_t)
	\end{cases}
\end{equation}
fails to admit a solution.
Such issues cannot occur in sample-and hold mode, unless $f(\cdot, U_t)$ or $\sigma(\cdot, U_t)$ is dicontinuous.
It is fairly intuitive that a control Lyapunov function constructed for a deterministic system will partially preserve its properties if some noise were to be added, unless the magnitude of that noise is too large.
However, naturally, the existence of a control Lyapunov function generally speaking does not imply much about stabilizability of systems described by \eqref{eqn:sys-core}, since the noise term $\sigma(X_t, U_t)Z_t\diff t$ can just entirely overhwelm the deterministic dynamics of the system. A way to remedy this problems is to consider bounded noise models.

Now, consider a special case of \eqref{eqn:sys-core}:
\begin{equation}
\label{eqn:sys-main}
	\begin{cases}
			& \diff X_t = f(X_t, U_t) \diff t + \sigma (X_t, U_t) Z_t \diff t, \ X_0 = x_0, \\ 
	& f \ : \ \mathbb{X} \times \mathbb{U} \rightarrow \mathbb{R}^n, \ \sigma \ : \ \mathbb{X} \times \mathbb{U} \rightarrow \mathbb{R}^{n \times d}, \\
	& Z_t \text{ is a measurable random process w.r.t. }t, \\
	& \forall t \in \mathbb{R} \spc \ \nrm{Z_t} \leq \Zb, \text{  (bounded noise)}\\
	& U_t := \mu(X_{t - (t \ \text{mod} \ \delta)}), \text{    (sample and hold policy)}
	\end{cases}
\end{equation}
where $\diff t$ implies Lebesgue integration.

We obtained the above system from \eqref{eqn:sys-core} by assuming bounded noise and asserting that a sample-and-hold policy is in place.


On the first glance it may seem like $\nrm{Z_t} \leq \Zb$ imposes a significant restriction on physical objects that this model can represent, since even Gaussian white noise fails to be described by it. 
However, upon a closer inspection it becomes clear that this theoretical loss of generality does not have much of a negative impact in practice. 
Note that noise represented by the following model can very closely mimic Brownian noise:
\begin{equation}
  \begin{aligned}
   & Z_t \diff t =  \xi_{\ceil{\lambda t}} \diff t, \\
   & \quad \text{where } \ \xi_i \sim \bar{\mathcal{N}}(0, \sigma^2, \Zb), \{\xi_i\} \text{ are independent}
    \end{aligned}
\end{equation}

Indeed, for real-world objects some values of noise are large enough to be considered impossible as opposed to merely improbable.
Therefore, an adequate reduction of support may in fact make the noise model more accurate. 

There are also a number of common noise models that imply this kind of boundedness \cite{Domingo2020-bounded-stochastic-processes}, in particular:
\begin{itemize}
	\item The Doering-Cai-Lin (DCL) noise \\
	\begin{equation}
	\label{eqn:DCL}
		\diff Z_t = - \tfrac 1 \theta Z_t \diff t + \sqrt{\tfrac {1 - Z_t^2 }{\theta (\gamma + 1 )}} \diff B_t,
	\end{equation}
	with parameters $\gamma>-1, \theta>0$;
	\item The Tsallis-Stariolo-Borland (TSB) noise \\
	\begin{equation}
		\label{eqn:TSB}
		\diff Z_t = - \tfrac 1 \theta \tfrac{Z_t}{1 - Z_t^2} \diff t + \sqrt{\tfrac{1 - q}{\theta}} \diff B_t,
	\end{equation}
	with $\theta > 0, \spc q < 1$ parameters;
	\item Kessler–S{\o}rensen (KS) noise \\
	\begin{equation}
		\label{eqn:KS}
		\diff Z_t = - \tfrac {\vartheta} {\pi \theta} \tan \left ( {\tfrac \pi 2 Z_t} \right ) \diff t + \tfrac {2} {\pi \sqrt{\theta (\gamma + 1 )}}  \diff B_t,
	\end{equation}   
	with $\theta > 0, \gamma \ge 0, \vartheta = \tfrac{2\gamma+1}{\gamma+1}$ parameters. 	
\end{itemize}

The following theorem demonstrates how a control Lyapunov function can be used to stabilize \eqref{eqn:sys-main} in sample-and-hold mode.
The theorem explicitly provides a stabilizing policy $\mu(x)$ for a given control Lyapunov function $L(x)$.
The theorem also provides estimates for the quality of resulting attraction.

\begin{thm}
\label{thm_bounded_noise}
	Let the following assumptions hold for system \eqref{eqn:sys-main}:
	\begin{framed}
	\begin{itemize}
		\item[(A1)] Bounds as per \eqref{eqn:core-bounds} exist.
		\item[(A2)] There exists a control Lyapunov function $L(\cdot)$:
	\begin{equation}
		\al{
		&\forall x \in \X \spc \inf\limits_{v \in \co(f(x, \U))} \D_v L(x) \le - \alpha_3(\nrm{x}),	\\
				&\forall x\in \mathbb{X} \ \ \alpha_1(\nrm{x}) \leq L(x) \leq \alpha_2(\nrm{x}),\\
		&\alpha_1, \alpha_2, \alpha_3 \in \K_\infty, \ L \spc \ \mathbb{X} \rightarrow \mathbb{R}_+.\\
		}
	\end{equation}

	\item[(A3)]
	The policy $\mu(x)$ is chosen so as to satisfy
	\begin{equation}
		\al{
		& \scal{\frac{x - {}^\lambda x}{\lambda^2},f({}^\lambda x, \mu(x))} \leq \\
		& \inf\limits_{u \in \mathbb{U}} \ \scal{\frac{x - {}^\lambda x}{\lambda^2},f({}^\lambda x, u)} + \eta, \\ 
		& \text{where ${}^\lambda x := \argmin\limits_{x^* \in \mathbb{X}}(2\lambda^2L(x^*) + \nrm{x^* - x})$}, \ \lambda > 0.
		}
	\end{equation}
	\item[(A4)]
	\begin{equation}
		\al{
		& \exists \ R > 0 \ \exists i \in \N \spc \ \\
		& \quad \quad \quad \quad y \left( \circ^{i}\tilde{r}_{\bar{Z}, \lambda, \delta, \eta}(R^* ); \delta\right) \leq R.
		}
	\end{equation}
	where $R^* := \alpha_1^{-1}(\alpha_2(R))$.
	\end{itemize}
	\end{framed}
	Then a unique global Caratheodory solution 
	almost surely exists and the following can be claimed:
	\begin{framed}
		\begin{itemize}
			\item[(C1)]
			$\lim\limits_{l\rightarrow \infty}r_l = \rst$ exists,
			\begin{equation}
				\text{where } r_l := \tilde{r}(r_{l-1}, \bar{Z}, \lambda, \delta, \eta), \ r_0 = R.
			\end{equation}
		
			\item[(C2)]
			The origin of the system is semi-asymptotically stable in probability in $R$ until $r$,
			\begin{equation}
				\label{eqn:def-r}
				\quad \quad \quad  \quad \text{where } r := \aaoi{\aat{y(\rst; \delta)}}.
			\end{equation}
		
			
			
			
			
			\item[(C3)]
			The origin of the system is semi-asymptotically stable on average in $R$ until $\rb$,
			\begin{equation}
				\al{
				\text{where } \rb := {}^\lambda \breve{\alpha}^{-1}_1({}^\lambda \invbreve{\alpha}_2(& \rht(r, \evb, \sdb, \lambda, \delta, \eta) +\\
				& + (\bar{f}(r) + \bar{\sigma}(r)\evb)\delta)), \\
				}\\
			\end{equation}
			 and ${}^\lambda \breve{\alpha}, {}^\lambda \invbreve{\alpha}_2$ are envelopes over [0, r].
			\end{itemize}
	\end{framed}
\end{thm}

\begin{crl}
	(C2) and (C3) would hold even if $\rst$ were to be replaced with any $r_l$.
	Therefore there is no need to compute $\rst$ precisely, since the theorem holds for any of its approximations $r_l$.
\end{crl}

\begin{rem}
	In the deterministic case, the analogous theorem would imply that by selecting a sufficiently small $\delta > 0$ we could make the attractor arbitrarily small and the basin of attraction arbitrarily large.
	All of the claims in (C1), (C2) and (C3) would hold with $\Zb$ substituted for $0$. 
	Furthermore, deterministic counterparts of this theorem can be found in \cite{Clarke1997-stabilization} and its generalization that accounted for computational uncertainty \cite{Osinenko2018practstabilization}.
	The related result under system and measurement uncertainty in deterministic form was further addressed in \cite{schmidt2021inf}.
	A particular difference to the latter work is that also convergence in mean in the sense of Definition \ref{dfn_conv_mean} is shown in Theorem \ref{thm_bounded_noise}.
\end{rem}

\begin{rem}
	The assumption (A3) means that $\mu(x)$ is chosen by minimizing $\scal{\frac{x - {}^\lambda x}{2\lambda^2},f({}^\lambda x, \mu(x))}$ and the minimization error is at most $\eta$. This kind of an approximate minimizer will always exist. 
	Furthermore the above theorem accounts for the fact that in practice we almost never have the luxury of using an exact minimizer.
\end{rem}

\begin{rem}
	It may seem that (A4) imposes restrictions on the systems, to which the theorem is applicable, however it can be derived that for sufficiently small $\bar{Z}$, $\delta$ and $\eta$ this assumption will hold.
	In other words, if (A4) does not hold, then that essentially means that noise magnitude, sampling time and optimization errors are too large to infer stability in $\ball_R$. 
\end{rem}

\begin{crl}
	For arbitrary $0 < r < R$, system \eqref{eqn:sys-main} is semi-asymptotically stable in probability in $R$ until $r$, provided that $\Zb$, $\delta$ and $\eta$ are sufficiently small.
\end{crl}

\begin{crl}
	For arbitrary $0 < r < R$, system \eqref{eqn:sys-main} is semi-asymptotically stable on average in $R$ until $r$, provided that $\E{\nrm{Z_t}^2}$, $\delta$ and $\eta$ are sufficiently small.
\end{crl}

Once an attraction function is obtained for a system of the kind \eqref{eqn:sys-main}, it can be used to inspect the influence of noise, latency and optimization error on stability. This implies a number of possible applications
\begin{itemize}
	\item \textbf{Verification.} Known values of noise, latency and optimization error can be used to infer (average) semi-asymptotic stability for given radii, producing a formal guarantee.
	Computing $r = \aaoi{\aat{y(r_l; \delta)}}$ for some $l$ will yield a radius in which the state will remain permanently once it gets there. 
	Likewise, one could compute $\rb = {}^\lambda \breve{\alpha}^{-1}_1({}^\lambda \invbreve{\alpha}_2(\rht(r, \evb, \sdb, \lambda, \delta, \eta) + (\bar{f}(r) + \bar{\sigma}(r)\evb)\delta))$ to determine a substantially smaller radius in which the mean state will stay.
	\item \textbf{Tuning.} The attraction function can be used to discover values of parameters that ensure desired quality of stabilization.
	If we were to  consider $r$ and $\rb$ as functions of $(\Zb, \delta, \eta)$ and $(\evb, \sdb, \delta, \eta)$ respectively, then both of those functions would be strictly increasing with respect to each one of their arguments.
	Thus suitable values of parameters can be identified through a simple grid search.
	\item \textbf{Analysis of robustness.} Given a Lyapunov function with a corresponding policy for a deterministic control system, one could investigate its robustness by inspecting how introducing noise, approximation error and response latency affects the quality of stabilization.
	This way, evaluating $r$ and $\rb$ for various sets of parameters will show in which settings this controller performs sufficiently well.
	The latter will reveal suitable implementations and environments, for which the system is guaranteed to function as intended.
	\item \textbf{Safe reinforcement learning.} If semi-asymptotic stability has been inferred over $R$, then any controller that implements $\mu(x)$ when near $R$ is guaranteed to stay within $\aaoi{\aat{R}}$ regardless of the policy used in other areas of the state space.
	This allows to utilize numerical optimization of reinforcement learning, while maintaining formal guarantees of Lyapunov theory (see \eg \cite{Osinenko2019stabactorcritic}).
\end{itemize}


\begin{rem}
	Note that there are no extra restrictions on $f(x, \cdot)$, $\sigma(x, \cdot)$ and $\mathbb{U}$, other than the ones imposed by (A1). 
	For the above theorem to hold $f(x, \cdot)$, $\sigma(x, \cdot)$ and $\mathbb{U}$ do not even have to be measurable. 
	By using approximate minimizers, we circumvent having to rely on the extreme value theorem.
\end{rem}


In the next section, we briefly discuss the case when the system trajectory is an It\^o process driven by a Brownian motion hence unbounded noise.

\section{It\^o Processes}
\label{sec_unbounded}

The case of unbounded noise can too be considered if assumptions of a different kind were to be made. Consider the following It\^o drift-diffusion process, driven by standard Brownian motion $B_t$: 
\begin{equation}
\label{eqn:sys-ito}
	\begin{cases}
	& \diff X_t = f(X_t, U_t)\diff t + \sigma(X_t, U_t)\diff B_t, \\
	& U_t := \mu(X_{t - (t \ \text{mod} \ \delta)}), \text{    (sample and hold policy)} \\
	\end{cases}
\end{equation} 
The generator of the stochastic differential equation of \eqref{eqn:sys-ito}, for a smooth function $L$, is defined as follows:
\begin{equation}
	\label{eqn:stoch-generator}
	\begin{aligned}
		\A^\mu L(x) = &\nabla L f(x, \mu(x)) + \\
		& + \frac 1 2 \tr{(\sigma(x, \mu(x)))^\top \nabla^2 L(x) \sigma(x, \mu(x))},		
	\end{aligned}
\end{equation}
where $\nabla L$ is the gradient vector and $\nabla^2 L$ is the Hessian \ie the matrix of second-order derivatives.
We define the operator $\Gamma^\mu_R$ as follows:
\begin{equation}
	\al{
	&\Gamma^\mu_R L := \bar f^\mu_{R} \big( \lip{R}{\nabla L} \bar f^\mu_{R} + \\ 
	& + b_{R} \Lip_f(R) + \bar \sigma_{R} b'_{R} \Lip_\sigma(R) + \frac 1 2 (\bar \sigma_{R}^{\mu})^2 \lip{\nabla^2 L}{R} \big).
	}
\end{equation}
where
\begin{equation}
	\al{
	& \bar f^\mu_R := \sup_{x \in \ball_R} \nrm{f(x, \mu(x))}, & \bar \sigma^\mu_R := \sup_{x \in \ball_R} \nrm{\sigma(x, \mu(x))},\\ 
	& b_R := \sup_{x \in \ball_R} \nrm{\nabla L(x)}, & b'_R := \sup_{x \in \ball_R} \nrm{\nabla^2 L(x)}.
	}
\end{equation}

Now that we can no longer rely on the boundedness of noise, to assert stability we have to demand the existence of a stronger version of the classical Lyapunov function.
\begin{dfn}[Stochastic Lyapunov pair]
	\label{dfn:stoch-L-pair}
	A stochastic Lyapunov pair $(L, \mu)$ is for the system $\eqref{eqn:sys-ito}$ a pair of functions if:
	$L \in \mathcal C^2$;
	there exist $\bar \alpha_1, \bar \alpha_2 > 0$, $\alpha_3 \in \K_\infty$ with $\alpha_3(\sqrt{\nrm{x}})$ convex; 
	there exists $\alpha_4 \in \K_\infty \spc \forall R \spc \Gamma^\mu_R L \le \alpha_4(R)$ \sut $\alpha_4(\sqrt{R})$ is concave;	
	(monotone condition cf. \cite{mao2007stochastic,yu2018sampled,fu2016sampled}) there exists $K > 0, K_{\mu} \spc \sut \forall x, \mu(x) \in \ball_{K_\mu}$ and $\forall x, u \in \ball_{K_\mu} \spc x^\top f(x, u) + \tfrac 1 2 \nrm{\sigma(x, u)}^2 \le K (1 + \nrm{x}^2)$; 
	the following properties hold:
	\begin{align}
		\label{eqn:LF-alphas12}
		& \forall x \spc \bar \alpha_1 \nrm{x}^2 \leq L(x) \leq \bar \alpha_2 \nrm{x}^2, \\
		\label{eqn:LF-alpha3}
		& \forall x \spc \A^\mu L(x) \le - \alpha_3(\nrm{x}) + \bar \Sigma, \bar \Sigma > 0.
	\end{align}
\end{dfn}

\begin{rem}
	The number $\bar \Sigma > 0$ in Definition \ref{dfn:stoch-L-pair} is related to the noise and differentiates the stochastic case from a deterministic one, which typically possesses a decay condition of the kind $\scal{\nabla L(x), f(x, \mu(x))} \le - \alpha_3(\nrm{x})$.
\end{rem}

\begin{rem}
	In general, for a function $\varphi$, the Jensen's gap \ie $\E{\varphi(X)} - \varphi( \E{X} )$ can be arbitrarily large.
	To relate various expected values in the analysis of Theorem \ref{thm_stoch-stab}, the Jensen's inequality has to be utilized, which motivates the stated conditions.  
	Notice \eg \cite[Lemma~2.1]{Reif1999-EKF-stab} also used quadratic bounding functions of $L$.
	A condition $\A^\mu L(x) = - c L + \bar \Sigma, \ c > 0$ (cf. \cite[Theorem~4.1]{Deng2001-stochastic-stab-noise}) also fits the assumptions since $- c L \le - c \bar \alpha_1 \nrm{x}^2$ and $\alpha_3$ is thus effectively $c \bar \alpha_1 \nrm{x}^2$ and so $\alpha_3(\sqrt{\nrm{x}})$ is convex.
	The function $\alpha_4(\sqrt{\nrm{x}})$ being concave is satisfied if \eg $L$ is quadratic, $f(\cdot, \mu(\cdot))$ is Lipschitz and of linear growth, $\sigma(\cdot, \mu(\cdot))$ is Lipschitz and bounded.
	This condition may be seen as restrictive, although \eg linear growth is often assumed in stochastic stabilization in mean (see, for instance, \cite{li2021stochastic,lan2017global,yang2011mean,fu2016sampled,yu2018sampled}).
	The monotone condition is in the style of \cite{mao2007stochastic} and is weaker than in related works on sampled stochastic stabilization (see \eg \cite{yu2018sampled,fu2016sampled}), whereas it should be noted that universal formulas with bounded controls are known \cite{Lin1991-stabilization-bounded-controls}.
	Furthermore, this condition will secure global existence of strong solutions \cite{mao2007stochastic}, which is unavoidable in case of \SH mode, and this is in contrast to the ``standard'' Lyapunov techniques in stochastic systems \cite{khasminskii201-stochastic}.
	The reason is that, in the latter, decay of the subject Lyapunov function is ensured for all times, whereas in the herein considered case, there are necessarily time intervals in which the said decay cannot be guaranteed.
%
%
%
%
\end{rem}

\begin{thm}[It\^o process semi-asymp. stab. on average]
	\label{thm_stoch-stab}
	Consider a stochastic system \eqref{eqn:sys-ito}.
	Suppose there exists a stochastic Lyapunov pair $(L, \mu)$ and $U_t = \mu(X_{t - (t\text{ mod }\delta)})$.
	Then for each $R > 0$, $\bar{\alpha}_3 > 0$  there exists a sufficiently small $\delta > 0$, that \eqref{eqn:sys-ito} is semi-asymptotically stable on average in $R$ until $\rho$, where 
	\begin{equation}
	\rho = \inf_{\rho'}\{ \inf_{r' \ge  \sqrt{\tfrac{\bar \alpha_1}{2 \bar \alpha_2}}\rho'  } (\alpha_3(r') - \bar \Sigma) \ge \bar \alpha_3 \}. 
	\end{equation}
\end{thm}

\begin{IEEEproof}
	See \cite[Theorem~1]{Osinenko2021stochSH}.
\end{IEEEproof}

\begin{crl}
	\label{crl:determ2stoch}
	Suppose that \eqref{eqn:LF-alpha3} has the form
	\begin{equation}
		\label{eqn:smooth-control-Lyapunov-function-decay-noiseless}
		\forall x \spc \nabla L(x) f(x,\mu(x)) \le - \alpha_3(\nrm{x}).
	\end{equation}
	In other words, $(L, \mu)$ is a Lyapunov pair for the noiseless system $\dot x = f(x, u)$.
	If it holds that
	\begin{equation}
		\label{eqn:sigma2nabla2-growth}
		\begin{aligned}
			& \exists \tilde r > 0 \spc \forall x \notin \ball_{\tilde r} \\
			& \alpha_3(\nrm{x}) \ge \tfrac 1 2 \nrm{\sigma^\top(x,\mu(x)) \nabla^2 L(x) \sigma(x,\mu(x))},
		\end{aligned}
	\end{equation}
	then assuming $U_t := \mu(X_{t - (t\text{ mod }\delta)})$, for each $R$ there exists a sufficiently small $\delta > 0$, that \eqref{eqn:sys-ito} is semi-asymptotically stable on average in $\rho$ over $R$
	
	In particular, \eqref{eqn:sigma2nabla2-growth} holds if $\nrm{\sigma}$ is uniformly bounded and $\nrm{\nabla^2 L(x)}$ has a growth rate lower than that of $\alpha_3$ everywhere except for a vicinity of the origin.
\end{crl}

\begin{rem}
	The growth condition \eqref{eqn:sigma2nabla2-growth} in Corollary \ref{crl:determ2stoch} may be justified as follows.
	Roughly speaking, taking derivatives decreases the growth rate.
	That is, one would normally expect that, outside some vicinity of the origin, $\nrm{\nabla^2 L(x)}$ grows slower than $\nrm{\nabla L(x)}$.
	Such is the case when $L$ is \eg polynomial.
	The diffusion function $\sigma$ describes the noise magnification depending on the state and control action.
	It may be justified in some applications to assume this term to be bounded uniformly in $x, u$.
	All in all, Corollary \ref{crl:determ2stoch} gives a particular hint on transferring a Lyapunov pair from a nominal, noiseless, to a noisy system.
\end{rem}

\section{Experimental study.}
\label{sec_experiment}

\begin{figure}[h]
	\vspace{2\parskip}
	\includegraphics[width=0.85\columnwidth]{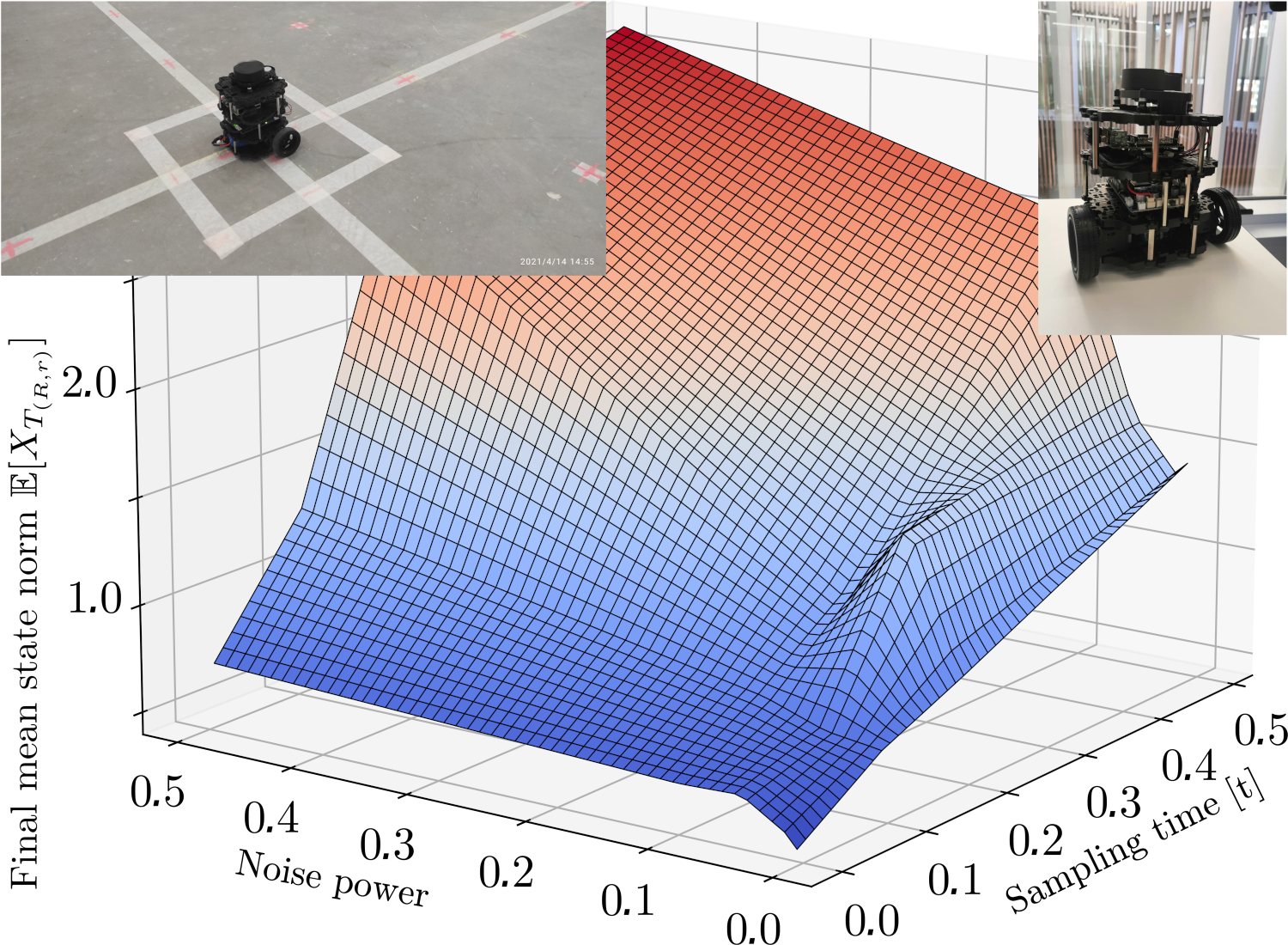}
	\caption{Effects of sampling time and noise power on practical stabilization of a mobile robot.
	Top left: fragment of the test bed.
	Top right: test robot TurtleBot 3.}
	\label{fig:experim}
\end{figure}

An experimental study of mobile robot parking was performed to demonstrate the effects of the developed theory (see Fig. \ref{fig:experim}).
Experiments were performed under practical truncated white Gaussian noise of varying power introduced into the system.

We used the above described method of inf-convolutions and the model $(\dot x_1, \dot x_2, \dot x_3)^\top = (u_1, u_2, x_1 u_2 - x_2 u_1 )^\top $ with a corresponding non-smooth control Lyapunov function $L(x) = x_1^2 + x_2^2 +2 x_3^2 - 2 |x_3|\sqrt{x_1^2 + x_2^2}$ \cite{Braun2017-feedback}, which is global and the control actions are confined to $[-1,1]^2$.
Notice that for a fixed noise power parameter, the decrease of the sampling time has limited effect due to the noise.

\begin{rem}
Theorem \ref{thm_bounded_noise} bears a fairly general character as physical systems are most adequately described by stochastic differential equations and controllers are commonly implemented in a sampled mode.
Further examples of such settings include \eg mechanical and robotic systems driven by digital micro-controllers, stock market bots with time-discrete decisions to buy or sell, medical devices like insulin pumps with scheduled administration, greenhouse with digital climate control etc.
The nature of noise in such examples can be interpreted in various ways \eg as unmodelled mechanical disturbance, say, due to friction or wear; noisy sensors; random fluctuation of stock prices; fluctuations in homeostasis; fluctuations in plant respiration etc.

\end{rem}




\section{Conclusion}
This work was concerned with stabilization of nonlinear dynamical systems in the sample-and-hold framework.
A novel theoretical result was derived that enables synthesis, verification, tuning and robustness testing of digital stabilizers for stochastic systems with bounded noise.

\appendix

\subsection{\textbf{Preliminaries}}
\subsubsection{\textbf{Existence of a radial growth bound}}

\begin{lem}
	\label{lem:rgb}
	Under assumption (A1) there exists a radial growth bound $\rho(\cdot)$.
\end{lem}

\begin{IEEEproof}
	Observe that
	\begin{equation}
		\al{
		&x^Tf(x, u) + \nrm{\sigma^T(x, u)x}\Zb \leq \\ 
		& \nrm{x}\nrm{f(x, u)} + \nrm{\sigma(x, u)}\nrm{x}\Zb \leq \\
		& \nrm{x}(\fb(\nrm{x}) + \sgb(\nrm{x})\Zb).}
	\end{equation}
	Let $\rho^-(\nrm{x}) := \fb(\nrm{x}) + \sgb(\nrm{x})\Zb$.
	Since the latter function is upper semi-continuous, by the extreme value theorem it has a maximum over any compact set.
	Now, we construct a locally Lipschitz continuous function $\rho(\cdot)$ that bounds $\rho^-(\cdot)$ from above as follows:
	\begin{equation}
	\al{
	\rho(y) = &(1 - y \ \text{mod} \ 1)\max\limits_{\tilde{y}\in [\floor{y}, \floor{y} + 1]}\rho^-(\tilde{y}) \ + \\ 
	& \ (y \ \text{mod} \ 1)\max\limits_{\tilde{y}\in [\floor{y} + 1, \floor{y} + 2]}\rho^-(\tilde{y}).
	}
	\end{equation}
	It is evident that $\rho(\cdot)$ is non-negative, locally Lipscitz continuous and $\rho(y) \geq \rho^-(y)$.
	Thus, we have
	\begin{equation}
	x^Tf(x, u) + \nrm{\sigma^T(x, u)x}\Zb \leq \nrm{x} \rho(\nrm{x}).
	\end{equation}
\end{IEEEproof}

\subsubsection{\textbf{Local existence of radial forecast}}

\begin{lem}
	\label{lem:rf}
	For every radial growth bound, there exists a corresponding radial forecast.
	\end{lem}
	
\begin{IEEEproof}
	According to Definition \ref{dfn:rgb}, any radial growth bound is locally Lipschitz continuous.
	Therefore by the Picard-Lindel\"{o}f theorem, (\ref{eqn:radial-forecast}) has a local solution for any non-negative $y(0):=y_0$.
	By Definition \ref{dfn:rf}, the latter local solution qualifies as a radial forecast $y(y_0; \cdot)$.
\end{IEEEproof}

\begin{lem}
	\label{lem:af-maf}
	(A1) implies the existence of an attraction function.
\end{lem}

\begin{IEEEproof}
	This follows from Lemma \ref{lem:rgb} and Lemma \ref{lem:rf}, since the attraction function is explicitly constructed from a radial forecast in \eqref{eqn:attraction-function}.
\end{IEEEproof}

\subsubsection{\textbf{Bounding properties of radial forecast}}

\begin{lem}
	Let $t_1 > 0$, then the inequality $\nrm{X_{t_1}} \leq R$ implies that
	\begin{equation}
	\forall t_2 \geq t_1 \spc \ t_2 \in \dom \ (y(R; \cdot)) \ \implies \ \nrm{X_{t_2}} \leq y(R; t_2),
	\end{equation}
	where $\dom$ denotes the domain of a function.
\end{lem}

\begin{IEEEproof}
	Observe that
	\begin{equation}
		\al{
		& \nrm{X_t} \diff \nrm{X_t} = X_t \cdot \diff X_t = \\ & (X_t^Tf(X_t, U_t) + X_t^T\sigma(X_t, U_t)Z_t) \diff t \leq \\ & \rho(\nrm{X_t})\nrm{X_t}\diff t.}
	\end{equation}
	Thus, 
	\begin{equation} 
		\label{eqn:identity-inequality}
		\al{
		& \diff \nrm{X_t} \leq \rho(\nrm{X_t})\diff t, \\
		& \diff y(t) = \rho(y(t)) \diff t,\\
		& X_{t_1} \leq y(t_1).
		}
	\end{equation}
	Let $h(t) := \nrm{X_{t}} - y(t)$. Now assume that $h(t_2) > 0$, then
	\begin{equation}
		\exists c\in[t_1, t_2) \spc \ h(c) = 0.
	\end{equation}
	Now, let $d := \sup \{a \in [t_1, t_2] \ | \ h(a) \leq 0\}$.
	Obviously, $d \geq c$ and $h(d) = 0$.
	By subtracting $\diff y(t) = \rho(y(t)) \diff t$ from $\diff \nrm{X_t} \leq \rho(\nrm{X_t})\diff t$ we get
	\begin{equation}
		\label{eqn:pre-gronwall}
		\al{
		\forall t\in[d, t_2] \spc \ & h(t) \leq h(d) + \int_d^t\rho(\nrm{X_\tau}) - \rho(y(\tau))\diff \tau =  \\ & \int_d^t\rho(\nrm{X_\tau}) - \rho(y(\tau))\diff \tau.}
	\end{equation}
	Now, let $\bar{X} := \sup\limits_{t \in [t_1, t_2]}\nrm{X_t}$.
	Note that since $y(R; \cdot)$ is non-decreasing and $X(t_2) > y(t_2)$, we have 
	\begin{equation}
		\forall t\in[t_1, t_2] \ y(\tau) \leq \bar{X}.
	\end{equation} 
	Therefore,
	\begin{equation}
		\al{
		& \forall t\in[d, t_2] \ \int_d^t\rho(\nrm{X_\tau}) - \rho(y(\tau))\diff \tau \leq \\ 
		& \int_d^t\Lip_\rho(\bar{X})\lvert h(\tau)\rvert \diff \tau = \int_d^t\Lip_\rho(\bar{X}) h(\tau) \diff \tau.}
	\end{equation}
	Combining this with \eqref{eqn:pre-gronwall} gives
	\begin{equation}
		\al{
		& \forall t\in[d, t_2] \ h(t) \leq \int_d^t\Lip_\rho(\bar{X}) h(\tau) \diff \tau.}
	\end{equation}
	Thus, by the Gr\"{o}nwall inequality we have
	\begin{equation}
		h(t_2) \leq 0,
	\end{equation}
	which contradicts our assumption that $h(t_2) > 0$
\end{IEEEproof}

\subsubsection{\textbf{Domain of definition remarks for attraction function}}

\begin{lem}
	Domain of definition of an attraction function occupies the entirety of $\mathbb{R}_+$ as $\delta$ approaches $0$.
\end{lem}

\begin{IEEEproof}
	Due to existence of local solutions
	\begin{equation}
		\forall R \geq 0 \ \exists \delta > 0 \spc \ \exists y(R; \delta)
	\end{equation}
	Thus 
	\begin{equation}
		\al{
		\forall R \geq 0 \ \exists \delta > 0 \spc \  & R \in \dom \ \tilde{r}(\cdot, \Zb, \lambda, \delta, \eta).}
	\end{equation}
\end{IEEEproof}

\subsection{\textbf{Decay}}
\subsubsection{\textbf{Infimal convolution bounds}}

\begin{lem}
	\label{lem:inf-conv-bound}
	The following inequality holds:
	\begin{equation}
		\nrm{{}^\lambda x - x}^2 \leq 2\lambda^2L(x).
	\end{equation}
\end{lem}

\begin{IEEEproof}
	Assume that $\nrm{{}^\lambda x - x}^2 > 2\lambda^2L(x)$, then
	\begin{equation}
	\La{x} > \La{{}^\lambda x} + L(x) > L(x),
	\end{equation}
	which is a contradiction.
\end{IEEEproof}

\subsubsection{\textbf{$L_\lambda$ Lyapunov properties}}

\begin{lem}
	The function $L_\lambda(\cdot)$ is positive-definite.
\end{lem}

\begin{IEEEproof}
	Obviously, $L_\lambda(0) = 0$. Now assume that if
	\begin{equation}
		\exists x\neq 0 \spc \ \La{x} \leq 0,
	\end{equation}
	then, 
	\begin{equation}
		L({}^\lambda x) + \frac{\nrm{{}^\lambda x - x}^2}{2\lambda^2} \leq 0.
	\end{equation}
	Since both terms in the left-hand side of the last displayed formula are non-negative, we have
	\begin{equation}
		\al{
		& L({}^\lambda x) = 0,\\
		& \frac{\nrm{{}^\lambda x - x}^2}{2\lambda^2} = 0,
		}
	\end{equation}
	which implies
	\begin{equation}
		{}^\lambda x = 0 \implies \frac{\nrm{x}^2}{2\lambda^2} = 0 \implies x=0,
	\end{equation}
	which in turn contradicts $x \neq 0$.
\end{IEEEproof}

\begin{lem}
	The function $\La{\cdot}$ is radially unbounded.
\end{lem}

\begin{IEEEproof}
	Suppose there is an unbounded sublevel set
	\begin{equation}
		A=\{x \ | \ \La{x} \leq c\}.
	\end{equation}
	Let $x_l \in A$ be an unbounded sequence.
	Then, note that
	\begin{equation}
		\al{
		 & \left\{^\lambda x_l\right\}_l \text{ is unbounded} \implies \left\{L(^\lambda x_l)\right\}_l \text{ is unbounded}, \\
		 & \left\{^\lambda x_l\right\}_l \text{ is bounded} \implies \left\{\frac{\nrm{ ^\lambda x_l - x_l}}{2\lambda^2}\right\}_l \text{ is unbounded},
		}
	\end{equation}
	whence $\La{x_l} = L(^\lambda x_l) + \frac{\nrm{^\lambda x_l - x_l}}{2\lambda^2}$ is an unbounded sequence, which contradicts $\La{x_l} \leq c$.
\end{IEEEproof}

\subsubsection{\textbf{Upper bound for Decay} $\dprod{\zeta_\lambda}{f({}^\lambda x, \mu(x))}$ }

We define a proximal subgradient of $L_\lambda$ at ${}^\lambda x$ as $\zeta_\lambda = \frac{x - {}^\lambda x}{2\lambda^2}$.

\begin{lem}
\label{lem:final-bound}
	From (A3) it follows that
	\begin{equation}
		\dprod{\zeta_\lambda}{f({}^\lambda x, \mu(x))} \leq -\alpha_3(\max(0, \nrm{x} -\lambda\sqrt{2\at{\nrm{x}}})) + \eta.
	\end{equation}
\end{lem}

\begin{IEEEproof}
	According to (A3):
	\begin{equation}
		\dprod{\zeta_\lambda }{f({}^\lambda x, \mu(x))} \leq \inf\limits_{u \in \mathbb{U}} \ \dprod{\zeta_\lambda }{f({}^\lambda x, u)} + \eta.
	\end{equation}
	It is known from \cite{Osinenko2018practstabilization} that
	\begin{equation}
		\dprod{\zeta_\lambda }{f({}^\lambda x, u)} \leq \mathcal{D}_{f({}^\lambda x, u)}L({}^\lambda x).
	\end{equation}
	Thus, 
	\begin{equation}
		\al{
		& \inf\limits_{u \in \mathbb{U}} \ \dprod{\zeta_\lambda }{f({}^\lambda x, u)} + \eta \leq \inf\limits_{u \in \mathbb{U}} \ \mathcal{D}_{f({}^\lambda x, u)}L({}^\lambda x) + \eta \leq \\ 
		& -\alpha_3(\nrm{{}^\lambda x}) + \eta.}
	\end{equation}
	By \textit{Lemma \ref{lem:inf-conv-bound}}, we have
	\begin{equation}
		-\alpha_3(\nrm{{}^\lambda x}) + \eta \leq  -\alpha_3(\nrm{x} -\lambda\sqrt{2\at{\nrm{x}}}) + \eta.
	\end{equation}
\end{IEEEproof}

Let there be some number $t' \in (0, \delta]$.
Now, let $F_k := X_{k\delta + t'} - X_{k\delta}$.

Notice that
\begin{equation}
\label{eqn:decay-bound-2}
	\al{
	& F_k = \intstepg{\rhs}  \implies \\
	& \nrm{F_k} \leq \intstepg{\nrm{\rhs}} \leq \\ &\delta(\fb\left(\max\limits_{\tau \in [k\delta, k\delta + t']}X_\tau\right) + \sgb\left(\max\limits_{\tau \in [k\delta, k\delta + t']}X_\tau\right)\Zb) \leq \\ 
	& \delta(\fb\left(y(\nrm{X_{k\delta}}; t')\right) + \sgb\left(y(\nrm{X_{k\delta}}; t')\right)\Zb).\\
	}
\end{equation}

From this point onward, we denote $\fb\left(y(\nrm{X_{k\delta}}; t')\right)$ and $\sgb\left(y(\nrm{X_{k\delta}}; t')\right)$ as simply $\fb$ and $\sgb$ accordingly.
Denote
\begin{equation}
\label{eqn:x-large-inf-conv}
	{}^\lambda X_t := \argmin\limits_{y \in \mathbb{X}}\left(L(y) + \frac{\nrm{y - X_t}^2}{2\lambda^2}\right).
\end{equation}

Let us determine an upper bound for $\dprod{\zeta_\lambda }{F_k}$.
Notice that $F_k$ can be expressed in the following way:
\begin{multline}
	\al{
	F_k  = & t' f(X_{k\delta}, U_{k\delta}) + \underbrace{\intstepg{f(X_t, U_{k\delta}) - f(X_{k\delta}, U_{k\delta})}}_{=:A_1} + \\ 
	& + \underbrace{\intstepg{\sigma(X_t, U_{k\delta})Z_t}}_{=:A_2}.}
\end{multline}

Using this expression, deduce:
\begin{equation}
\label{eqn:bound-dot-product-1}
	\al{
	& \dprod{\zeta_\lambda }{F_k} = t'\dprod{\zeta_\lambda }{f(X_{k\delta}, U_{k\delta})} + \dprod{\zeta_\lambda }{A_1 + A_2}. \\}
\end{equation}

For the latter term in the right-hand side of the last displayed relation, we have
\begin{equation}
	\al{
	& \dprod{\zeta_\lambda }{A_1 + A_2} \leq \nrm{\zeta_\lambda }(\nrm{A_1} + \nrm{A_2}),\\
	& \nrm{A_1} \leq \Lip_f(y(\nrm{X_{k\delta}}; t'))(\fb + \sgb\Zb)\frac{{t'}^2}{2},\\
	& \nrm{A_2} \leq \sgb\Zb t'.\\
	}
\end{equation}
whereas for the former term, we have
\begin{equation}
\label{eqn:decay-term-bound-1}
	\al{
	& \dprod{\zeta_\lambda }{f(X_{k\delta}, U_{k\delta})} = \dprod{\zeta_\lambda }{f(X_{k\delta}, U_{k\delta}) - f(^\lambda X_{k\delta}, U_{k\delta})} + \\
	& + \dprod{\zeta_\lambda }{f(^\lambda X_{k\delta}, U_{k\delta})}, \\
	& \dprod{\zeta_\lambda }{f(X_{k\delta}, U_{k\delta}) - f(^\lambda X_{k\delta}, U_{k\delta})} \leq \\ 
	& \Lip_f\left(\nrm{X_{k\delta}} +\lambda\sqrt{2\at{\nrm{X_{k\delta}}}}\right)\nrm{\zeta_\lambda }\nrm{X_{k\delta} - {}^\lambda X_{k\delta}}.\\
	}
\end{equation}

Using \eqref{eqn:x-large-inf-conv}, the factor $\nrm{\zeta_\lambda }\nrm{X_{k\delta} - {}^\lambda X_{k\delta}}$ can be expressed and bounded as follows:
\begin{equation}
\label{eqn:decay-term-bound-1-1}
	\al{
	&\nrm{\zeta_\lambda }\nrm{X_{k\delta} - {}^\lambda X_{k\delta}} = 2(\La{X_{k\delta}} - L({}^\lambda X_{k\delta})) \leq \\ 
	& 2(L(X_{k\delta}) - L({}^\lambda X_{k\delta})) \leq \\ 
	& 2\Lip_L\left(\nrm{X_{k\delta}} +\lambda\sqrt{2\at{\nrm{X_{k\delta}}}}\right)\nrm{X_{k\delta} - {}^\lambda X_{k\delta}} \leq \\
	& 2\Lip_L\left(\nrm{X_{k\delta}} +\lambda\sqrt{2\at{\nrm{X_{k\delta}}}}\right)\lambda\sqrt{2\at{\nrm{X_{k\delta}}}}.
	}
\end{equation}

Substituting the obtained bounds, we get
\begin{equation}
\label{eqn:decay-bound-1}
	\al{
	&\dprod{\zeta_\lambda }{F_k} \leq 
	t'\Big(\dprod{\zeta_\lambda }{f(^\lambda X_{k\delta}, U_{k\delta})} + \\ 
	& 2\Lip_L\big(w(\nrm{X_{k\delta}})\big)\Lip_f\big(w(\nrm{X_{k\delta}})\big)\lambda\sqrt{2\at{\nrm{X_{k\delta}}}}\Big)  \\ 
	& + \frac{\sqrt{2\alpha_2(\nrm{X_{k\delta}})}}{\lambda}\Big(\frac{{t'}^2}{2}\Lip_f(y(\nrm{X_{k\delta}}; t'))(\fb + \sgb\Zb) \\ 
	& + t'\sgb\Zb\Big),
	}
\end{equation}
where $w(\nrm{X_{k\delta}}) = \nrm{X_{k\delta}} +\lambda\sqrt{2\at{\nrm{X_{k\delta}}}}$.

It is known from \cite{Osinenko2018practstabilization} that
\begin{equation}
\label{eqn:taylor-nonsmooth}
	\La{X_{k\delta + t'}} - \La{X_{k\delta}} \leq \dprod{\zeta_\lambda }{F_k} + \frac{\nrm{F_k}^2}{2\lambda^2}.
\end{equation}

Using (\ref{eqn:decay-bound-2}) and (\ref{eqn:decay-bound-1}) we get
\begin{equation} 
\label{eqn:decay-final}
	\al{
	&\La{X_{k\delta + t'}} - \La{X_{k\delta}} \leq \\
	& -t'\alpha_3\big(\max(0, \nrm{X_{k\delta}} -\lambda\sqrt{2\at{\nrm{X_{k\delta}}}})\big) + t'\nu(\nrm{X_{k\delta}}, t'),
	}
\end{equation}
where
\begin{equation} 
	\al{
	&\nu(\nrm{X_{k\delta}}, t') = \eta + \\ 
	& + \frac{\sqrt{2\alpha_2(\nrm{X_{k\delta}})}}{\lambda}\Big(2\Lip_L\big(w(\nrm{X_{k\delta}})\big)\Lip_f\big(w(\nrm{X_{k\delta}})\big)\lambda^2 + \\
	&+ \frac{t'}{2}\Lip_f(y(\nrm{X_{k\delta}}; t'))(\fb + \sgb\Zb) + \sgb\Zb\Big)  + \frac{t'(\fb + \sgb\Zb)^2}{2\lambda^2}.
	}
\end{equation}
Note that the attraction function is by definition
\begin{equation}
\label{eqn:attraction-short}
	\rtd(R, \Zb, \lambda, \delta, \eta) = \adi{\nu(R, \delta)} +\lambda\sqrt{2\at{R}}.
\end{equation}

\subsubsection{\textbf{Radial decay}}

Let $x_l$ denote $X_{l\delta}$.

\begin{lem}
\label{lem:3.1.1}
	If
	\begin{equation}
		\rtd := \rtd(R, \Zb, \lambda, \delta, \eta) < R,
	\end{equation}
	then, the following condition holds:
	\begin{equation}
		\al{
		&\forall \eps>0 \ \exists\Delta(\eps)<0 \ \forall k \in \mathbb{Z}  \ \forall t' \in (0, \delta] \spc \\
		& \nrm{X_{k\delta}} \in [\rtd + \eps, R] \implies \La{X_{k\delta + t'}} - \La{X_{k\delta}} \leq t'\Delta(\eps).}
	\end{equation}
\end{lem}

\begin{IEEEproof}
	Observe that
	\begin{equation}
		\al{
		& \nrm{X_{k\delta}} \geq \rtd(R, \Zb, \lambda, \delta, \eta) + \eps \geq \rtd(R, \Zb, \lambda, t', \eta) + \eps\implies \\
		& \ad{\nrm{X_{k\delta}} - \eps -\lambda\sqrt{2\at{R}}} \geq \\
		& \ad{\rtd(R, \Zb, \lambda, t', \eta) - \lambda\sqrt{2\at{R}}} \geq  \nu(R, t') \geq \\
		& \geq \nu(\nrm{X_{k\delta}}, t'). \\
		}
	\end{equation}
	Therefore, using \eqref{eqn:decay-final} we get
	\begin{equation}
		\al{
		& \ad{\nrm{X_{k\delta}} - \eps -\lambda\sqrt{2\at{R}}} - \ad{\nrm{X_{k\delta}} -\lambda\sqrt{2\at{R}}} \geq \\
		& \nu(\nrm{X_{k\delta}}, t') - \ad{\nrm{X_{k\delta}} -\lambda\sqrt{2\at{\nrm{X_{k\delta}}}}} \geq \\
		& \frac{\La{\nrm{X_{k\delta + t'}}} - \La{\nrm{X_{k\delta}}}}{t'}.}
	\end{equation}
	If we denote $z := \nrm{X_{k\delta}} - \eps -\lambda\sqrt{2\at{R}}$, we obtain
	\begin{equation}
	\label{eqn:core-decay}
		\al{
		& \frac{\La{\nrm{X_{k\delta + t'}}} - \La{\nrm{X_{k\delta}}}}{t'} \leq \\
		& \sup\limits_{z \in [\rtd - \tsqrt, R - \eps - \tsqrt]} (\ad{z} - \ad{z + \eps}) \leq \\ 
		& -\inf\limits_{z\in[0, R]}(\ad{z + \eps} - \ad{z}) < 0,
		}
	\end{equation}
	since $\ad{\cdot}$ is strictly increasing. We conclude the proof by simply asserting $\Delta(\eps) := -\inf\limits_{z\in[0, R]}(\ad{z + \eps} - \ad{z})$, which gives us
	\begin{equation}
		\La{\nrm{X_{k\delta + t'}}} - \La{\nrm{X_{k\delta}}} \leq t'\Delta(\eps).
	\end{equation}
\end{IEEEproof}

Let us denote $\rtd_{ \Zb, \lambda, \delta, \eta}(R)$ as simply $\rtd(R)$.

\subsubsection{\textbf{Iterated decay}}

\begin{lem}
	\label{lem:3.2.1}
	If we assume that $\rtd(R) < R$, then
	\begin{equation}
		\forall l \in \mathbb{N} \ \exists \bar{\eps} > 0 \ \forall \eps \in (0, \bar{\eps}) \spc \ \rtd(\circ^{l}\rtd(R) + \eps) < \circ^{l}\rtd(R) + \eps.
	\end{equation}
\end{lem}

\begin{IEEEproof}
	Since $\adi{\nu(\cdot, \delta)}$ is composed of non-decreasing functions and $\lambda\sqrt{2\at{\cdot}}$ is strictly increasing, due to \eqref{eqn:attraction-function} attraction function $\rtd(\cdot)$ is strictly increasing.
	Therefore, $\rtd(R) < R$ implies
	$\forall l \in \mathbb{N} \spc \ \circ^{l}\rtd(R) < \circ^{l - 1}\rtd(R)$.
	Now, obviously
	\begin{equation}
		\al{
		& \forall 0<\eps<\circ^{l - 1}\rtd(R) - \circ^{l}\rtd(R) \spc \\
		& \rtd(\circ^{l}\rtd(R) + \eps) < \rtd(\circ^{l - 1}\rtd(R)) < \circ^{l}\rtd(R) + \eps.
		}
	\end{equation}
\end{IEEEproof}

Attraction function $\rtd(\cdot)$ is composed of upper semi-continuous and non-decreasing continuous functions, therefore it is upper semi-continuous.
Thus $\rtd(\cdot)$ is right continuous, since it is both strictly increasing and upper semi-continuous.

\begin{lem}
\label{lem:rcont-implication}
	\begin{equation}
		\al{
		&\exists \bar{\eps}_3 > 0 \ \forall \eps_3\in(0, \bar{\eps}_3) \ \exists \eps_1\in (0, \bar{\eps}_1) \ \exists \eps_2\in (0, \bar{\eps}_2) \spc \\
		& \circ^{l + 1}\rtd(R) + \eps_3 = \rtd(\circ^{l}\rtd(R) + \eps_1) + \eps_2.
		}
	\end{equation}
\end{lem}

\begin{IEEEproof}
	This follows from the right-continuity of $\rtd(\cdot)$.
\end{IEEEproof}

\begin{lem}
	The condition $\rtd(R) < R$ implies
	\begin{equation}
		\al{
		& \forall l\in\mathbb{N} \ \forall 0 < \eps < R-\circ^{l}\rtd(R) \ \exists \Delta^l(\eps)<0 \spc \\ 
		&  \nrm{x_k} \in [\circ^{l}\rtd(R) + \eps, R] \ \implies \\
		& \La{x_{k + 1}} - \La{x_k} \leq \delta\Delta_l(\eps). 
		}
	\end{equation}
\end{lem}

\begin{IEEEproof}
	The above statement has already been proven for $n = 1$, so let us assume that it holds for some $n$:
	\begin{equation}
		\al{
		S_l \iff \ & \forall 0 < \eps_1 < R-\circ^{l}\rtd(R) \ \exists \Delta_l(\eps_1)<0 \spc \\ 
		&  \nrm{x_k} \in [\circ^{l}\rtd(R) + \eps_1, R] \ \implies \\
		& \La{x_{k + 1}} - \La{x_k} \leq \delta\Delta_l(\eps_1).}
	\end{equation}
	By \textit{Lemma \ref{lem:3.2.1}}, we have $\rtd(\circ^{l}\rtd(R) + \eps_1) < \circ^{l}\rtd(R) + \eps_1$, thus by \textit{Lemma \ref{lem:3.1.1}},
	\begin{equation}
		\al{&\forall 0<\eps_2<\circ^{l}\rtd(R) + \eps_1 - \rtd(\circ^{l}\rtd(R) + \eps_1) =: \bar{\eps}_2 \ \exists \Delta(\eps_2) < 0  \spc \\
		& \nrm{x_k} \in [\rtd(\circ^{l}\rtd(R) +\eps_1)+\eps_2, \circ^{l}\rtd(R) + \eps_1] \implies \\
		& \La{x_{k + 1}} - \La{x_k} \leq \delta\Delta(\eps_2).}
	\end{equation}
	By \textit{Lemma \ref{lem:rcont-implication}}, this implies
	\begin{equation}
		\al{
		& \forall \eps_3\in(0, \bar{\eps}_3) \ \exists \Delta(\eps_3)<0 \spc \\
		& \nrm{x_k} \in [\circ^{l + 1}\rtd(R) + \eps_3, \circ^{l}\rtd(R) + \eps_1],
		}
	\end{equation}
	which in conjunction with $S_l$ yields
	\begin{equation}
		\al{
		& \forall \eps_3\in(0, \bar{\eps}_3) \ \exists \tilde{\Delta}_{l + 1}(\eps_3)<0 \spc \\
		& \nrm{x_k} \in [\circ^{l + 1}\rtd(R) + \eps_3, R] \implies \\
		& \La{x_{k + 1}} - \La{x_k} \leq \delta\tilde{\Delta}_{l + 1}(\eps_3),
		}
	\end{equation}
	where $\tilde{\Delta}_{l + 1}(\eps_3) := \min(\Delta(\eps_3), \Delta_l(\eps_1(\eps_3)))$.
	
	Consider $\Delta_{l + 1}(\eps_3) := \tilde{\Delta}_{l + 1}(\min(\eps_3, \frac{\bar{\eps}_3}{2}))$.
	Then,
	\begin{equation}
		\al{
		S_{l + 1} \iff \ & \forall 0 < \eps_3 < R-\circ^{l + 1}\rtd(R) \ \exists \Delta_{l + 1}(\eps_1)<0 \spc \\ 
		&  \nrm{x_k} \in [\circ^{l + 1}\rtd(R) + \eps_3, R] \ \implies \\
		& \La{x_{k + 1}} - \La{x_k} \leq \delta\Delta_{l + 1}(\eps_3).}
	\end{equation}
	Thus, the induction step is concluded.
\end{IEEEproof}

\subsubsection{\textbf{Limit decay}}

Since $\circ^{l}\rtd(R)$ is a decreasing, it has a limit.
Let us denote this limit as $\rst$.

\begin{lem}
\label{lem:3.3.1}
	If $\rtd(R) < R$ then
	\begin{equation}
		\al{
		& \forall \eps \in(0, R - \rst) \ \exists \Delta^*(\eps) < 0 \spc \\
		& \nrm{x_k} \in [\rst + \eps, R] \implies \\
		& \La{x_{k + 1}} - \La{x_k} \leq \delta\Delta^*(\eps).
		}
	\end{equation}
\end{lem}

\begin{IEEEproof}
	The existence of a limit implies
	\begin{equation}
		\forall \eps \in(0, R - \rst) \ \exists l\in \mathbb{N} \ \circ^{l}\rtd(R) < \rst + \eps.
	\end{equation}
	Thus,
	\begin{equation}
		\al{
		& \circ^{l}\rtd(R) + \frac{\rst + \eps - \circ^{l(\eps)}\rtd(R)}{2} < \rst + \eps \implies \\
		& \nrm{x_k} \in [\rst + \eps, R] \implies \\ 
		& \La{x_{k + 1}} - \La{x_k} \leq \delta\Delta_{l(\eps)}\left(\frac{\rst + \eps - \circ^{l(\eps)}\rtd(R)}{2}\right).
		}
	\end{equation}
\end{IEEEproof}

\begin{lem}
\label{lem:decreasing}
	An attraction function $\rtd(\cdot)$ has the following property:
	\begin{equation}
		\forall r\in(\rst, R] \spc \ \rtd(r) < r.
	\end{equation}
\end{lem}

\begin{IEEEproof}
	The statement is already proven for $r = R$.
	Assume that 
	\begin{equation}
		\exists r^+ \in (\rst, R) \spc \ \rtd(r^+) \geq r^+,
	\end{equation}
	but due to convergence
	\begin{equation}
		\exists l \in \mathbb{N} \spc \ (\circ^{l}\rtd(R) > r^+) \text{ and } (\circ^{l + 1}\rtd(R) < r^+),
	\end{equation}
	which implies
	\begin{equation}
		\rtd(\circ^{l}\rtd(R)) > \rtd(r^+) \geq r^+ > \circ^{l + 1}\rtd(R),
	\end{equation}
	an obvious contradiction.
\end{IEEEproof}

\begin{crl}
	The number $\rst$ is the largest root of $\rtd(r) - r$ over $[0, R]$.
\end{crl}

\subsection{\textbf{Semi-asymptotic stability in probability}}
\subsubsection{\textbf{Attraction}}

Let $\Rst$ denote $\aaoi{\aat{R}}$.
Obviously $R \leq \Rst$.
\begin{lem}
	\label{lem:4.1.1}
	Under the assumption that
	\begin{equation}
		\exists i \in \mathbb{N} \spc \ \circ^{i}\rtd(\Rst) < R,
	\end{equation}
	we have
	\begin{equation}
		\lim\limits_{l\rightarrow \infty}\circ^{l}\rtd(\Rst) = \lim\limits_{l\rightarrow \infty}\circ^{l}\rtd(R).
	\end{equation}
\end{lem}

\begin{IEEEproof}
	Notice $\circ^{l}\rtd(\Rst)$ is a bounded decreasing sequence, and, therefore, it has a limit.
	Furthermore, it has a subsequence $\circ^{l}\rtd(\circ^{i}\rtd(\Rst))$ that is bounded above by $\circ^{l}\rtd(R)$, while the sequence $\circ^{l}\rtd(\Rst)$ itself is bounded below by $\circ^{l}\rtd(R)$.
	Therefore,
	\begin{equation}
	\rst \leq \lim\limits_{l\rightarrow \infty}\circ^{l}\rtd(\Rst) \leq \rst.
	\end{equation}
\end{IEEEproof}

\begin{lem}
	\label{lem:4.1.2}
	If we assume that
	\begin{equation}
	\exists i \in \mathbb{N} \spc \ \circ^{i}\rtd(\Rst) < R,
	\end{equation}
	then, the following holds
	\begin{equation}
	\nrm{x_k} \leq R \rightarrow \forall \eps > 0 \ \exists l \geq k \spc \ \nrm{x_l} < \rst + \eps.
	\end{equation}
\end{lem}

\begin{IEEEproof}
	First note that $\rtd(\Rst) < \Rst$, because the negation of this statement implies $\circ^{i + 1}\rtd(\Rst) \geq \circ^{i}\rtd(\Rst)$, which contradicts \textit{Lemma \ref{lem:decreasing}}.
	
	Now let us assume the opposite of the statement is to be proven:
	\begin{equation}
		(\nrm{x_k} \leq R) \text{ and } \exists \eps > 0 \ \forall l \geq k \spc \ \nrm{x_l} \geq \rst + \eps.
	\end{equation}
	
	Since $\rtd(\Rst) < \Rst$, \textit{Lemma \ref{lem:3.3.1}} and \textit{Lemma \ref{lem:4.1.1}} yield
	\begin{equation}
		\al{
		& \forall \eps \in(0, R - \rst) \ \exists \Delta^*(\eps) < 0 \spc \\
		& \nrm{x_k} \in [\rst + \eps, \Rst] \implies \\
		& \La{x_{k + 1}} - \La{x_k} \leq \delta\Delta^*(\eps).
		}
	\end{equation}
	Thus, if we assume that $\La{x_l} \leq \aat{R}$, then 
	\begin{equation}
	\al{
	& \La{x_{l + 1}} \leq \La{x_l} + \Delta^*(\eps) \leq \aat{R} \implies \\
	& \nrm{x_{l + 1}} \leq \Rst,
	}
	\end{equation}
	which in turn proves $\forall l\geq k \spc \ \nrm{x_l} \leq \Rst$ by induction, which then implies
	\begin{equation}
	\La{x_l} \leq \La{x_k} + (n - k)\delta\Delta^*(\eps).
	\end{equation}
	Naturally, for $n > k - \frac{\La{x_k}}{\delta\Delta^*(\eps)}$ this yields $\La{x_l} < 0$, which is a contradiction.
\end{IEEEproof}

\subsubsection{\textbf{Ultimate boundedness}}

\begin{lem}
	\label{lem:4.2.1}
	Under the assumption that
	\begin{equation}
		\al{
		& \exists i \in \mathbb{N} \spc \ \circ^{i}\rtd(\Rst) < R, \\
		& \rst < y(\rst; \delta) \leq R, \\
		& \nrm{x_k} \leq R. \\
		}
	\end{equation}
	the following holds
	\begin{equation}
		\exists N \geq k \ \forall l \geq N \ \nrm{x_l} \leq \aaoi{\aat{y(\rst; \delta)}}=\bar{r}.
	\end{equation}
\end{lem}

\begin{IEEEproof}
	Note, that
	\begin{equation}
		y(\rst; \delta) = \rst + (y(\rst; \delta) - \rst) = \rst + \eps.
	\end{equation}
	By \textit{Lemma \ref{lem:4.1.2}} $\exists l \in \mathbb{N}\cup\{0\}  \ \nrm{x_l} < \rst + \eps$.
	Now let us assume that for some $m \geq n$
	\begin{equation}
		\La{x_i} \leq \aat{y(\rst; \delta)}.
	\end{equation}
	Then, 
	\begin{equation}
		\al{
		& \nrm{x_i} \leq \rst \implies \nrm{x_{i + 1}} \leq y(\rst; \delta) \implies \\
		& \implies \La{x_{i + 1}} \leq \aat{y(\rst; \delta)},
		}
	\end{equation}
	whereas 
	\begin{equation}
		\al{
		&\nrm{x_i} > \rst \implies \La{x_{i + 1}} \leq \La{x_{i}} \leq \aat{y(\rst; \delta)}.
		}
	\end{equation}
	Since the statement holds for $m = n$, we have
	\begin{equation}
		\al{
		&\forall m\geq l(\eps) \spc \ \La{x_{i}} \leq \aat{y(\rst; \delta)} \implies \\
		&\implies \forall m\geq l(\eps) \spc \ \nrm{x_i} \leq \aaoi{\aat{y(\rst; \delta)}}.
		}
	\end{equation}
\end{IEEEproof}

\begin{crl}
	If $y(\rst; \delta) = \rst$, then $\limsup\limits_{l\rightarrow \infty}\nrm{x_l} \leq \rst$.
\end{crl}

\begin{lem}
\label{lem:almost-sure-ultimate-boundedness}
\label{lem:C2}
	Assume that 
	\begin{equation}
		\al{
		& \exists i \in \mathbb{N} \spc \ \circ^{i}\rtd(\Rst) < R, \\
		& \rst < y(\rst; \delta) \leq R, \\
		& \nrm{x_k} \leq R. \\
		}
	\end{equation}
	Then the following holds
	\begin{equation}
		\exists t_2 \geq k\delta \ \forall t \geq t_2 \spc \ \nrm{X_t} \leq \aaoi{\aat{y(\rst; \delta)}}.
	\end{equation}
\end{lem}

\begin{IEEEproof}
	By \textit{Lemma \ref{lem:4.2.1}} 
	\begin{equation}
		\exists N \geq k \ \forall l \geq N \ \nrm{x_l} \leq \aaoi{\aat{y(\rst; \delta)}},
	\end{equation}
	Now consider an arbitrary $t \geq N\delta$, and let $i := t - (t \text{ mod } \delta)$. Then, using Lemma \ref{lem:3.3.1} we get
	\begin{equation}
		\al{
		& \nrm{x_{i}} > \rst \implies \La{X_t} \leq \La{x_i}, \\
		& \nrm{x_{i}} \leq \rst \implies \nrm{X_t} \leq y(\rst; t \text{ mod } \delta).
		}
	\end{equation}
\end{IEEEproof}

\begin{crl}
	If $y(\rst; \delta) = \rst$, then $\limsup\limits_{t\rightarrow \infty}\nrm{X_t} \leq \rst$.
\end{crl}

\begin{lem} 
	The norm state $\nrm{X_t}$ is of bounded variation on an arbitrary segment $[a, b]$.
\end{lem}

\begin{IEEEproof}
	Note, that
	\begin{equation}
		\nrm{X_b - X_a} \leq (\fb + \sgb\Zb)|b - a|.
	\end{equation}
	therefore $X_t$ is locally Lipschitz continuous, which ensures that it has a finite total variation on any segment.
\end{IEEEproof}

\subsubsection{\textbf{Mean decay}}

Let $\Var{Z_t} \leq \sdb$ and $\E{\nrm{Z_t}} \leq \evb$.
\begin{lem}
	If we assume that for some $r$ 
	\begin{equation}
	\label{eqn:invariance-assumption}
		\forall t' \in (0, \delta] \spc \nm{X_{k\delta + t'}} \leq r,
	\end{equation}
	then, the following holds:
	\begin{equation}
		\forall t' \in (0, \delta] \spc \E{\nrm{F_k}^2} \leq 4((\fb(r) + \sgb(r)\evb)^2 + \sgb(r)^2\sdb^2){t'}^2.
	\end{equation}
\end{lem}

\begin{IEEEproof}
	Let $F(t) := X_t - k\delta$.  
	\begin{equation}
		\al{
		& \diff \nrm{F(t)} \leq \fb+\sgb\nrm{Z_t}\diff t \implies \\
		& \diff \nrm{F(t)}^2 \leq 2\nrm{F(t)}(\fb+\sgb\nrm{Z_t})\diff t \implies \\
		& \E{\nrm{F(t)}^2} \leq 2\int_{k\delta}^{t}\E{\nrm{F(\tau)}(\fb + \sgb\nrm{Z_\tau})}\diff \tau \implies \\
		& \E{\nrm{F(t)}^2} \leq 2\int_{k\delta}^{t}\sqrt{\E{\nrm{F(\tau)}^2}}\sqrt{\E{(\fb + \sgb\nrm{Z_\tau})^2}}\diff \tau.
		}
	\end{equation} 
	Let $\xi_t := \sqrt{\E{\nrm{F(t)}^2}}$, $\chi_t := \sqrt{\E{(\fb + \sgb\nrm{Z_\tau})^2}}$ then 
	\begin{equation}
		\xi^2_t \leq 2\int_{k\delta}^{t}\xi_\tau \chi_\tau\diff \tau.
	\end{equation}
	Now, let $\xi_{\max} := \sup\limits_{\tau\in[k\delta, t]}\xi_\tau$.
	Since $[k\delta, t]$ is compact and $\xi_{\tau}$ is continuous, we have $\xi_{\max} = \xi_{t_{\max}}$, where $t_{\max} \in [k\delta, t]$.
	This yields
	\begin{equation}
		\al{
		& \xi^2_{\max} \leq 2\int_{k\delta}^{t_{\max}}\xi_\tau \chi_\tau\diff \tau \implies \\
		& \xi^2_{\max} \leq 2\int_{k\delta}^{t_{\max}}\xi_{\max} \chi_\tau\diff \tau \implies \\ 
		& \xi_{\max} \leq 2\int_{k\delta}^{t_{\max}}\chi_\tau\diff \tau \implies \\ 
		& \xi_{\max} \leq 2\int_{k\delta}^{t}\chi_\tau\diff \tau \implies \\ 
		& \xi_t \leq 2\int_{k\delta}^{t}\chi_\tau\diff \tau \implies \\ 
		& \xi_t \leq 2\int_{k\delta}^{t}\sqrt{\fb(r) + 2\fb(r)\sgb\evb + \sgb(r)^2(\evb^2 + \sdb^2)}\diff \tau \implies \\ 
		& \xi^2_{k\delta + t'} \leq 4((\fb(r) + \sgb(r)\evb)^2 + \sgb(r)^2\sdb^2){t'}^2.
		}
	\end{equation}
\end{IEEEproof}

Now let us determine a bound for $\E{\dprod{\zeta_\lambda , F_k}}$ under the assumption \eqref{eqn:invariance-assumption}.
Recall from \eqref{eqn:bound-dot-product-1} that
\begin{equation}
	\al{
	& \dprod{\zeta_\lambda }{F_k} \leq \delta\dprod{\zeta_\lambda }{f(x_k, u_k)} + \nrm{\zeta_\lambda }(\nrm{A_1} + \nrm{A_2}) \leq \\
	& \delta\dprod{\zeta_\lambda }{f(x_k, u_k)} + \frac{\sqrt{2\at{\nrm{x_k}}}}{\lambda}(\nrm{A_1} + \nrm{A_2}). \\
	}
\end{equation}
Here, we consider the following bounds for $\nrm{A_1}$, $\nrm{A_2}$:
\begin{equation}
	\al{
	& \nrm{A_1} \leq \intstepg{\Lip_f\int_{k\delta}^t(\fb(r) + \sgb(r)\nrm{Z_\tau})\diff \tau},\\
	& \nrm{A_2} \leq \intstepg{\sgb(r)\nrm{Z_t}}.\\
	}
\end{equation}
This gives us
\begin{equation}
	\al{
	& \E{\nrm{A_1}} \leq \Lip_f(r)(\fb(r) + \sgb(r)\evb)\frac{{t'}^2}{2},\\
	& \E{\nrm{A_2}} \leq \sgb(r)\evb\delta.\\
	}
\end{equation}
Recall the bounds from \eqref{eqn:decay-term-bound-1} and \eqref{eqn:decay-term-bound-1-1}:
\begin{equation}
	\al{
	& \dprod{\zeta_\lambda }{f(x_k, u_k)} \leq \dprod{\zeta_\lambda }{f({}^\lambda x_k, u_k)} + \\ 
	& + \nrm{\zeta_\lambda }\nrm{f(x_k, u_k) - f({}^\lambda x_k, u_k)},\\
	& \nrm{\zeta_\lambda }\nrm{f(x_k, u_k) - f({}^\lambda x_k, u_k)} \leq \\
	& \nrm{\zeta_\lambda }\Lip_f\nrm{{}^\lambda x_k - x_k} \leq 2\Lip_f\Lip_L\lambda\sqrt{2\at{\nrm{x_k}}}.\\
	}
\end{equation}
We now apply $\E{\cdot}$ and obtain
\begin{equation}
	\al{
	& \E{\dprod{\zeta_\lambda }{F_k}} \leq t'\E{\dprod{\zeta_\lambda }{f({}^\lambda x_k, u_k)}} + \\ 
	& 2t' \lambda^2\Lip_f(r + \lambda \sqrt{2\at{r}})\Lip_L(r + \\
	& \lambda \sqrt{2\at{r}})\sqrt{2\at{\nrm{x_k}}} + \\ 
	& \frac{\sqrt{2\at{\nrm{x_k}}}}{\lambda}(\frac{{t'}^2}{2}\Lip_f(r)(\fb(r) + \sgb(r)\evb) + t'\sgb(r)\evb).
	}
\end{equation}
Recall that $\adc{\cdot}$ is defined as the lower convex envelope of $\ad{\cdot}$ over $[0, r]$.

\begin{lem}
\label{lem:final-bound-mean}
	If $(\E{\nrm{x_k}} \geq\lambda\sqrt{2\at{\nrm{x_k}}}) \text{ and } (\nrm{x_k} \leq r)$, then
	\begin{equation}
		\al{
		& \E{\dprod{\zeta_\lambda }{f({}^\lambda x_k, u_k)}} \leq \\
		& -\adc{\E{\nrm{x_k}} -\lambda\sqrt{2\at{r}}} + \eta.
		}
	\end{equation}
\end{lem}

\begin{IEEEproof}
	This follows from \textit{Lemma \ref{lem:final-bound}} in conjunction with the Jensen's inequality.
\end{IEEEproof}

From (\ref{eqn:taylor-nonsmooth}), we have
\begin{multline}
	\E{\La{X_{k\delta + t'}}} - \E{\La{X_{k\delta}}} \leq \\ 
	\E{\dprod{\zeta_\lambda }{F_k}} + \E{\frac{\nrm{F_k}^2}{2\lambda^2}}.
\end{multline}
Which, under the assumption that $(\E{\nrm{x_k}} \geq\lambda\sqrt{2\at{\nrm{x_k}}}) \text{ and } (\nrm{x_k} \leq r)$, expands to 
\begin{equation}
\label{eqn:mean-decay-final}
	\al{
	& \E{\La{x_{k + 1}}} - \E{\La{x_k}} \leq \\ 
	& -t'\adc{\E{\nrm{x_k}} -\lambda\sqrt{2\at{r}}} + \hat{\nu}(r, t')t', 
	}
\end{equation}
where
\begin{equation}
	\al{
	& \hat{\nu}(r, t') = \eta + \\
	& \sqrt{2\at{r}}\Big(2a\Lip_L(w(r))\Lip_f(w(r)) + \\
	& \frac{t'}{2a}\Lip_f(r)(\fb(r) + \sgb(r)\evb) + \\ 
	& \frac{\sgb(r)\evb}{\lambda}\Big) + \frac{2t'}{\lambda^2}((\fb(r) + \sgb(r)\Zb)^2 +\\
	& \sdb^2\sgb(r)^2).
	}
\end{equation}

\subsubsection{\textbf{Mean radial decay}}

The mean attraction function is defined as 
\begin{equation}
	\rht(r, \evb, \sdb, a, \delta, \eta) = \adci{\hat{\nu}(r, t')} +\lambda\sqrt{2\at{r}}.
\end{equation}

\begin{lem}
\label{lem:mean-decay}
	If $\rht(r) < r$, then
	\begin{equation}
		\al{
		& \forall \eps>0 \ \exists \Delta(\eps)< 0 \ \forall t' \in (0, \delta ] \spc \\
		& (\E{\nrm{x_k}} \geq \rht + \eps) \text{ and } (\nrm{x_k} \leq r) \implies \\
		& \E{\La{X_{k\delta + t'}}} - \E{\La{X_{k\delta}}} \leq t'\Delta(\eps).
		}
	\end{equation}
\end{lem}

\begin{IEEEproof}
	Note that $\E{\nrm{x_k}} \geq \rht + \eps$ implies 
	\begin{equation}
		\E{\nrm{x_k}} \geq\lambda\sqrt{2\at{r}} + \eps.
	\end{equation}
\end{IEEEproof}

Similarly to the proof of \textit{Lemma \ref{lem:3.1.1}}, we have
\begin{equation}
	\al{
	& \E{\nrm{x_k}} \geq \rht + \eps \implies \\
	& \implies \adc{\E{\nrm{x_k}} - \eps -\lambda\sqrt{2\at{r}}} \geq \hat{\nu}(r, t').\\
	}
\end{equation}
Using \eqref{eqn:mean-decay-final}, we get
\begin{equation}
	\al{
	& \adc{\nrm{x_k} - \eps -\lambda\sqrt{2\at{r}}} - \adc{\nrm{x_k} -\lambda\sqrt{2\at{r}}} \geq \\
	& \geq  \frac{\E{\La{X_{k\delta + t'}}} - \E{\La{X_{k\delta}}}}{t'}. \\
	}
\end{equation}
Finally, by analogy with \eqref{eqn:core-decay}, we obtain
\begin{equation}
	\al{
	& \adc{\nrm{x_k} - \eps -\lambda\sqrt{2\at{r}}} - \adc{\nrm{x_k} -\lambda\sqrt{2\at{r}}} \leq \\
	& -\inf\limits_{z\in [0, r]}(\adc{z + \eps} - \adc{z}) < 0,
	}
\end{equation}
which similarly allows us to assert $\Delta(\eps) := -\inf\limits_{z\in [0, \rb]}(\adc{z + \eps} - \adc{z})$.

\subsection{\textbf{Semi-asymptotic stability on average}}
\subsubsection{\textbf{Mean attraction}}

\begin{lem}
	If we assume that $\forall l \geq k \spc \ \nrm{x_l} \leq r$, then the following holds:
	\begin{equation}
		\forall \eps>0 \ \exists l\geq k \spc \ \E{\nrm{x_l}} < \rht(r) + \eps.
	\end{equation}
\end{lem}

\begin{IEEEproof}
	Assume that $\rht(r) \geq r$.
	Then, obviously $\E{\nrm{x_l}} < \rht(\rb) + \eps$.
	Thus let us assume that $\rht(r) \geq r$.
	Then, by \textit{Lemma \ref{lem:mean-decay}}, we have
	\begin{equation}
		\E{\nrm{x_l}} \geq \hat{r} + \eps \implies \E{\La{x_{l + 1}}} - \E{\La{x_l}} \leq \delta\Delta(\eps).
	\end{equation}
	But, if we assume the opposite of the statement to be proven, we have
	\begin{equation}
		\al{
		& \forall l\geq k \spc \ \E{\nrm{x_l}} \geq \hat{r} + \eps \implies \\
		& \forall l\geq k \spc \ \E{\La{x_{l + 1}}} - \E{\La{x_l}} \leq \Delta(\eps) \implies \\
		& \forall l > k - \frac{\E{\La{x_{k}}}}{\Delta(\eps)} \spc \ \E{\La{x_l}} < 0,
		}
	\end{equation}
	which is an obvious contradiction.
\end{IEEEproof}

\subsubsection{\textbf{Mean ultimate boundedness}}

Let $\aaoc{\cdot}$ be the upper concave envelope of $\aao{\cdot}$ over $[0, r]$ and let $\aatc{\cdot}$ be the lower convex envelope of $\aat{\cdot}$ over $[0, r]$.

\begin{lem}
\label{lem:mean-boundedness}
	If we assume that $\forall l \geq k \spc \ \nrm{x_l} \leq r$, then the following holds:
	\begin{equation}
		\al{
		& \exists N \geq k \ \forall l\geq N \spc \\
		& \E{\nrm{x_l}} \leq \aaoci{\aatc{\rht + (\fb(r) + \sgb(r)\evb)\delta}}.
		}
	\end{equation}
\end{lem}

\begin{IEEEproof}
	The condition $\aaoci{\aatc{\rht + \delta(\fb(r) + \sgb(r)\evb)}} = \rht$ would imply $\rb \leq \rht$, which immediately proves the lemma, so let us instead assume that $\aaoci{\aatc{\rht + (\fb(r) + \sgb(r)\evb)\delta}} > \rht$.
	This yields
	\begin{equation}
		\al{
		& \rht + (\fb(r) + \sgb(r)\evb)\delta = \rht + \eps \\ \implies
		& \exists N \geq k \spc \ \E{\nrm{x_l}} \leq \rht + (\fb(r) + \sgb(r)\evb)\delta.
		}
	\end{equation}
	By the Jensen's inequality, we have
	\begin{equation}
		\E{\La{x}} \leq \aatc{\E{\nrm{x}}}.
	\end{equation}
	Now, assume that for some $n \geq N$
	\begin{equation}
		\E{\La{x_l}} \leq \aatc{\rht + \delta(\fb(r) + \sgb(r)\cdot\evb)}.
	\end{equation}
	This means
	\begin{equation}
		\al{
		& \E{\nrm{x_l}} > \rht \implies \E{\La{x_{l + 1}}} < \E{\La{x_l}} \leq \\
		& \aatc{\rht + (\fb(r) + \sgb(r)\evb)\delta},
		}
	\end{equation}
	but at the same time
	\begin{equation}
		\al{
		& \E{\nrm{x_l}} \leq \rht \implies \E{\nrm{x_{l + 1}}} - \E{\nrm{x_l}} \leq \\
		& \E{\intstep{\fb(r) + \sgb(r)\nrm{Z_t}}} \implies \\
		& \E{\La{x_{l + 1}}} \leq \aatc{\E{\nrm{x_{l + 1}}}} \leq \\
		& \aatc{\rht + (\fb(r) + \sgb(r)\evb)\delta}.
		}
	\end{equation}
	Since the statement holds for $n = N$, the above constitutes a proof by induction of
	\begin{equation}
		\al{
		&\forall l \geq N \spc \ \E{\La{x_l}} \leq \aatc{\rht + \delta(\fb(r) + \sgb(r) \cdot\evb)}.
		}
	\end{equation}
	The statement of the lemma is evident if we consider the fact that
	\begin{equation}
		\aaoc{\E{\nrm{x_l}}} \leq \E{\La{x_l}}.
	\end{equation}
\end{IEEEproof}

\begin{lem}
\label{lem:C3}
	If we assume that
	\begin{equation}
		\al{
		& \exists i \in \mathbb{N} \spc \ \circ^{i}\rtd(\Rst) < R, \\
		& \rst < y(\rst; \delta) \leq R, \\
		& \nrm{x_k} \leq R. \\
		}
	\end{equation}
	Then, it follows that
	\begin{equation}
		\al{
		& \exists N \geq k \ \forall t\geq N\delta \spc \\
		& \E{\nrm{X_t}} \leq \aaoci{\aatc{\rht + \delta(\fb(r) + \sgb(r)\cdot\evb)}}.
		}
	\end{equation}
\end{lem}

\begin{IEEEproof}
	This follows from \textit{Lemma \ref{lem:C2}}, \textit{Lemma \ref{lem:mean-boundedness}} and \textit{Lemma \ref{lem:mean-decay}}.
\end{IEEEproof}

\subsection{\textbf{Proof of Theorem \ref{thm_bounded_noise}}}

\begin{itemize}
	\item[(C1)] From (A4) it follows that $\circ^{i + 1}\rtd(R) \leq \circ^{i}\rtd(R)$. 
	Since the sequence $\{\circ^{i}\rtd(R)\}_i$ is bounded from below, the sequence converges.
	\item[(C2)] Follows from (A4) by \textit{Lemma \ref{lem:C2}}.
	\item[(C3)] Follows from (A4) by \textit{Lemma \ref{lem:C3}}.
\end{itemize}

\bibliography{
bib/adversarial,
bib/analysis,
bib/automated-reasoning,
bib/computable,
bib/constr-math,
bib/ctrl-history,
bib/diff-games,
bib/discont-DE,
bib/Filippov-sol,
bib/formal-ctrl,
bib/lin-ctrl,
bib/logic,
bib/model-reduction,
bib/MPC,
bib/nonlin-ctrl,
bib/nonsmooth-analysis,
bib/opt-ctrl,
bib/Osinenko,
bib/perturb-thr,
bib/selector-misc,
bib/semiconcave,
bib/sensitivity,
bib/set-thr,
bib/sliding-mode,
bib/soft,
bib/stabilization,
bib/stochastic,
bib/topology,
bib/uncertainty,
bib/verified-integration,
bib/viability,
bib/metric-spaces,
bib/DP,
bib/constructing-LFs,
bib/Kalman-filter,
bib/RL
} 

\begin{thebibliography}{10}
\providecommand{\url}[1]{#1}
\csname url@samestyle\endcsname
\providecommand{\newblock}{\relax}
\providecommand{\bibinfo}[2]{#2}
\providecommand{\BIBentrySTDinterwordspacing}{\spaceskip=0pt\relax}
\providecommand{\BIBentryALTinterwordstretchfactor}{4}
\providecommand{\BIBentryALTinterwordspacing}{\spaceskip=\fontdimen2\font plus
\BIBentryALTinterwordstretchfactor\fontdimen3\font minus
  \fontdimen4\font\relax}
\providecommand{\BIBforeignlanguage}[2]{{%
\expandafter\ifx\csname l@#1\endcsname\relax
\typeout{** WARNING: IEEEtran.bst: No hyphenation pattern has been}%
\typeout{** loaded for the language `#1'. Using the pattern for}%
\typeout{** the default language instead.}%
\else
\language=\csname l@#1\endcsname
\fi
#2}}
\providecommand{\BIBdecl}{\relax}
\BIBdecl

\bibitem{khasminskii201-stochastic}
R.~Khasminskii and G.~Milstein, \emph{Stochastic Stability of Differential
  Equations}, ser. Stochastic Modelling and Applied Probability.\hskip 1em plus
  0.5em minus 0.4em\relax Springer, 2011.

\bibitem{Kushner1965stochastic-stability}
H.~J. Kushner, ``On the stability of stochastic dynamical systems,''
  \emph{Proc. National Academy of Sciences of USA}, vol.~53, no.~1, pp. 8--12,
  1965.

\bibitem{mao1991-stability}
X.~Mao, \emph{Stability of Stochastic Differential Equations with Respect to
  Semimartingales}, ser. Pitman research notes in mathematics series.\hskip 1em
  plus 0.5em minus 0.4em\relax Longman Scientific \& Technical, 1991.

\bibitem{Deng2001-stochastic-stab-noise}
H.~{Deng}, M.~{Krstic}, and R.~J. {Williams}, ``Stabilization of stochastic
  nonlinear systems driven by noise of unknown covariance,'' \emph{IEEE Tran.
  on Autom. Control}, vol.~46, no.~8, pp. 1237--1253, 2001.

\bibitem{Khalil1996-nonlin-sys}
H.~Khalil, \emph{Nonlinear {S}ystems}.\hskip 1em plus 0.5em minus 0.4em\relax
  Prentice-Hall. 2nd edition, 1996.

\bibitem{McAllister2003-stochastic}
R.~D. McAllister and J.~B. Rawlings, ``Stochastic {L}yapunov functions and
  asymptotic stability in probability,'' Tech. Rep., 08 2020.

\bibitem{Liu2008-stochastic}
S.-J. Liu, J.-f. Zhang, and Z.-p. Jiang, ``A notion of stochastic
  input-to-state stability and its application to stability of cascaded
  stochastic nonlinear systems,'' \emph{Acta Mathematicae Applicatae Sinica},
  vol.~24, pp. 141--156, 03 2008.

\bibitem{Liu2011-stochastic}
L.~Liu and X.-J. Xie, ``Output-feedback stabilization for stochastic high-order
  nonlinear systems with time-varying delay,'' \emph{Automatica}, vol.~47,
  no.~12, pp. 2772--2779, 2011.

\bibitem{Li2017-Stochastic}
H.~{Li}, L.~{Bai}, Q.~{Zhou}, R.~{Lu}, and L.~{Wang}, ``Adaptive fuzzy control
  of stochastic nonstrict-feedback nonlinear systems with input saturation,''
  \emph{IEEE Tran. on Syst., Man, and Cybernetics: Systems}, vol.~47, no.~8,
  pp. 2185--2197, 2017.

\bibitem{Wu2013-stochastic}
Z.~{Wu}, M.~{Cui}, P.~{Shi}, and H.~R. {Karimi}, ``Stability of stochastic
  nonlinear systems with state-dependent switching,'' \emph{IEEE Tran. on
  Autom. Control}, vol.~58, no.~8, pp. 1904--1918, 2013.

\bibitem{Cheng2001-stochastic-stability}
{Cheng Tang} and T.~{Basar}, ``Stochastic stability of singularly perturbed
  nonlinear systems,'' in \emph{IEEE CDC}, vol.~1, 2001, pp. 399--404 vol.1.

\bibitem{Teel2014-stochastic-stability}
A.~R. Teel, A.~Subbaraman, and A.~Sferlazza, ``Stability analysis for
  stochastic hybrid systems: A survey,'' \emph{Automatica}, vol.~50, no.~10,
  pp. 2435--2456, 2014.

\bibitem{Huang2008-hybrid-stochastic-retarded}
L.~Huang, X.~Mao, and F.~Deng, ``Stability of hybrid stochastic retarded
  systems,'' \emph{IEEE Tran. Circuits and Syst. I: Regular Papers}, vol.~55,
  no.~11, pp. 3413--3420, 2008.

\bibitem{Feng1992-stochastic-stability}
X.~{Feng}, K.~A. {Loparo}, Y.~{Ji}, and H.~J. {Chizeck}, ``Stochastic stability
  properties of jump linear systems,'' \emph{IEEE Tran. on Autom. Control},
  vol.~37, no.~1, pp. 38--53, 1992.

\bibitem{caraballo2015-practical-stochastic}
T.~Caraballo, M.~A. Hammami, and L.~Mchiri, ``On the practical global uniform
  asymptotic stability of stochastic differential equations,'' 2015.

\bibitem{Qin2020-stochastic}
Y.~Qin, M.~Cao, and B.~D.~O. Anderson, ``Lyapunov criterion for stochastic
  systems and its applications in distributed computation,'' \emph{IEEE Tran.
  on Autom. Control}, vol.~65, no.~2, pp. 546--560, 2020.

\bibitem{Do2020-stochastic}
K.~Do, ``Practical asymptotic stability of stochastic systems driven by {L}évy
  processes and its application to control of {TORA} systems,'' \emph{Int. J.
  Control}, pp. 1--25, 03 2020.

\bibitem{fu2016sampled}
X.~Fu, Y.~Kang, and P.~Li, ``Sampled-data stabilization for a class of
  stochastic nonlinear systems based on the approximate discrete-time models,''
  in \emph{2016 Australian Control Conference (AuCC)}.\hskip 1em plus 0.5em
  minus 0.4em\relax IEEE, 2016, pp. 258--263.

\bibitem{yu2018sampled}
P.~Yu, Y.~Kang, and Q.~Zhang, ``Sampled-data stabilization for a class of
  stochastic nonlinear systems with markovian switching based on the
  approximate discrete-time models,'' in \emph{2018 Australian \& New Zealand
  Control Conference (ANZCC)}.\hskip 1em plus 0.5em minus 0.4em\relax IEEE,
  2018, pp. 413--418.

\bibitem{GAO2018}
``Event-triggered control for stochastic nonlinear systems,''
  \emph{Automatica}, vol.~95, pp. 534--538, 2018.

\bibitem{yang2022stability}
X.~Yang and H.~Li, ``Stability analysis of probabilistic boolean networks with
  switching discrete probability distribution,'' \emph{IEEE Transactions on
  Automatic Control}, 2022.

\bibitem{li2020robustness}
H.~Li, X.~Yang, and S.~Wang, ``Robustness for stability and stabilization of
  boolean networks with stochastic function perturbations,'' \emph{IEEE
  Transactions on Automatic Control}, vol.~66, no.~3, pp. 1231--1237, 2020.

\bibitem{Brockett1983-stabilization}
R.~Brockett, ``Asymptotic stability and feedback stabilization,''
  \emph{Differential geometric control theory}, vol.~27, no.~1, pp. 181--191,
  1983.

\bibitem{Sontag1990-stabilization-survey}
E.~Sontag, ``Feedback stabilization of nonlinear systems,'' in \emph{Robust
  control of linear systems and nonlinear control}.\hskip 1em plus 0.5em minus
  0.4em\relax Springer, 1990, pp. 61--81.

\bibitem{Clarke1997-stabilization}
F.~Clarke, Y.~Ledyaev, E.~Sontag, and A.~Subbotin, ``Asymptotic controllability
  implies feedback stabilization,'' \emph{IEEE Tran. on Autom. Control},
  vol.~42, no.~10, pp. 1394--1407, 1997.

\bibitem{Fontes2003discontinuous}
F.~Fontes, ``Discontinuous feedbacks, discontinuous optimal controls, and
  continuous-time model predictive control,'' \emph{International Journal of
  Robust and Nonlinear Control: IFAC-Affiliated Journal}, vol.~13, no. 3-4, pp.
  191--209, 2003.

\bibitem{Cortes2008-discont-dyn-sys}
J.~Cortes, ``Discontinuous dynamical systems,'' \emph{IEEE {C}ontrol {S}ys.},
  vol.~28, no.~3, 2008.

\bibitem{Clarke2009-slid-mode-stab}
F.~Clarke and R.~Vinter, ``Stability analysis of sliding-mode feedback
  control,'' \emph{Control and Cybernetics}, vol.~4, no.~38, pp. 1169--1192,
  2009.

\bibitem{Braun2017-SH-stabilization-Dini-aim}
P.~Braun, L.~Gr{\"u}ne, and C.~Kellett, ``Feedback design using nonsmooth
  control {L}yapunov functions: A numerical case study for the nonholonomic
  integrator,'' in \emph{IEEE Conference on Decision and Control}, 2017.

\bibitem{Clarke2011-discont-stabilization}
F.~Clarke, ``{L}yapunov functions and discontinuous stabilizing feedback,''
  \emph{Annual Reviews in Control}, vol.~35, no.~1, pp. 13--33, 2011.

\bibitem{Baier2012-linear}
R.~Baier, L.~Gr{\"u}ne, and S.~Hafstein, ``Linear programming based {L}yapunov
  function computation for differential inclusions,'' \emph{Discrete and
  Continuous Dynamical Systems-Series B}, vol.~17, no.~1, pp. 33--56, 2012.

\bibitem{Baier2014-num-CLF}
R.~Baier and S.~Hafstein, ``Numerical computation of control {L}yapunov
  functions in the sense of generalized gradients,'' in \emph{Int. Symposium on
  Mathematical Theory of Networks and Systems}, 2014, pp. 1173--1180.

\bibitem{grammatico2013control}
S.~Grammatico, F.~Blanchini, and A.~Caiti, ``Control-sharing and merging
  control {L}yapunov functions,'' \emph{IEEE Transactions on Automatic
  Control}, vol.~59, no.~1, pp. 107--119, 2013.

\bibitem{Giesl2015-review}
P.~Giesl and S.~Hafstein, ``Review on computational methods for {L}yapunov
  functions,'' \emph{Discrete and Continuous Dynamical Systems-Series B},
  vol.~20, no.~8, pp. 2291--2331, 2015.

\bibitem{Osinenko2019nonsmoothstabsurvey}
P.~Osinenko, P.~Schmidt, and S.~Streif, ``Nonsmooth stabilization and its
  computational aspects,'' \emph{IFAC-PapersOnLine}, vol.~53, no.~2, pp.
  6370--6377, 2020, presented at IFAC World Congress.

\bibitem{Sontag1989-formula}
E.~Sontag, ``A "universal" construction of artstein's theorem on nonlinear
  stabilization,'' \emph{Syst. \& Control Letters}, vol.~13, no.~2, pp.
  117--123, 1989.

\bibitem{Lin1991-stabilization-bounded-controls}
Y.~Lin and E.~D. Sontag, ``A universal formula for stabilization with bounded
  controls,'' \emph{Syst. \& Control Letters}, vol.~16, no.~6, pp. 393--397,
  1991.

\bibitem{Florchinger1994-Lyapunov-like-stochastic-stability}
P.~{Florchinger}, ``Lyapunov-like techniques for stochastic stability,'' in
  \emph{IEEE CDC}, vol.~2, 1994, pp. 1145--1150.

\bibitem{Osinenko2021stochSH}
P.~Osinenko and G.~Yaremenko, ``On stochastic stabilization of sampled
  systems,'' in \emph{Conference on Decision and Control (CDC)}, 2021,
  accepted.

\bibitem{Domingo2020-bounded-stochastic-processes}
D.~Domingo, A.~d’Onofrio, and F.~Flandoli, ``Properties of bounded stochastic
  processes employed in biophysics,'' \emph{Stochastic Analysis and Appl.},
  vol.~38, no.~2, pp. 277--306, 2020.

\bibitem{Osinenko2018practstabilization}
P.~Osinenko, L.~Beckenbach, and S.~Streif, ``Practical sample-and-hold
  stabilization of nonlinear systems under approximate optimizers,'' \emph{IEEE
  Control Systems Letters}, vol.~2, no.~4, pp. 569--574, 2018, presented at
  Conference on Decision and Control (CDC).

\bibitem{schmidt2021inf}
P.~Schmidt, P.~Osinenko, and S.~Streif, ``On inf-convolution-based robust
  practical stabilization under computational uncertainty,'' \emph{IEEE
  Transactions on Automatic Control}, 2021.

\bibitem{Osinenko2019stabactorcritic}
P.~Osinenko, L.~Beckenbach, T.~G{\"o}hrt, and S.~Streif, ``A reinforcement
  learning method with real-time closed-loop stability guarantee,''
  \emph{IFAC-PapersOnLine}, vol.~53, no.~2, pp. 8043--8048, 2020, presented at
  IFAC World Congress.

\bibitem{mao2007stochastic}
X.~Mao, \emph{Stochastic differential equations and applications}.\hskip 1em
  plus 0.5em minus 0.4em\relax Elsevier, 2007.

\bibitem{Reif1999-EKF-stab}
K.~Reif, S.~Gunther, E.~Yaz, and R.~Unbehauen, ``{S}tochastic stability of the
  discrete-time extended {K}alman filter,'' \emph{IEEE Transactions on
  Automatic Control}, vol.~44, no.~4, pp. 714--728, 1999.

\bibitem{li2021stochastic}
W.~Li and M.~Krstic, ``Stochastic nonlinear prescribed-time stabilization and
  inverse optimality,'' \emph{IEEE Tran. on Autom. Control}, 2021.

\bibitem{lan2017global}
Q.~Lan and S.~Li, ``Global output-feedback stabilization for a class of
  stochastic nonlinear systems via sampled-data control,'' \emph{Int. J. Robust
  and Nonlin. Control}, vol.~27, no.~17, pp. 3643--3658, 2017.

\bibitem{yang2011mean}
S.~Yang, B.~Shi, and M.~Li, ``Mean square stability of impulsive stochastic
  differential systems,'' \emph{Int. J. Differential Equations}, vol. 2011,
  2011.

\bibitem{Braun2017-feedback}
P.~{Braun}, L.~{Grüne}, and C.~M. {Kellett}.

\end{thebibliography}

\bibliographystyle{IEEEtran}

\vfill\null
\pagebreak
\begin{IEEEbiography}
	[{\includegraphics[width=1in,height=1.25in,clip,keepaspectratio]{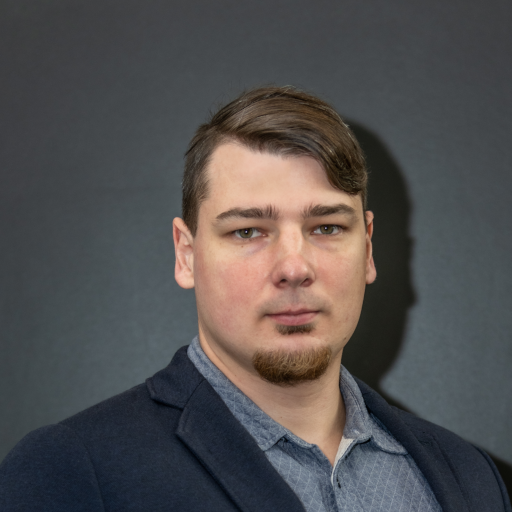}}]{Pavel Osinenko}
	studied control engineering at the Bauman Moscow State Technical University, Russia from 2003 through 2009.
	He started his postgraduate education in vehicle control systems at the same university.
	In 2011, Pavel moved to the Dresden University of Technology, Germany after receiving a Presidential Grant.
	He obtained a PhD degree in 2014 with a dissertation on vehicle optimal control and identification.
	Pavel has work experience in the German private sector and at the Fraunhofer Institute for Transportation and Infrastructure Systems in Dresden, Germany.
	In 2016, he made a transition to the Chemnitz University of Technology, Germany, as a postdoctoral researcher.
	He was responsible for project and research coordination, doctorate supervision, teaching, administration etc.
	Since 2020, Pavel has been an assistant professor at the Skolkovo Institute of Science and Technology, Russia, where he leads the Artificial Intelligence in Dynamic Action Lab.
	Pavel’s areas of research include reinforcement learning, especially its safety and connections to control theory, and computational aspects of dynamical systems.
\end{IEEEbiography}

\begin{IEEEbiography}
	[{\includegraphics[width=1in,height=1.25in,clip,keepaspectratio]{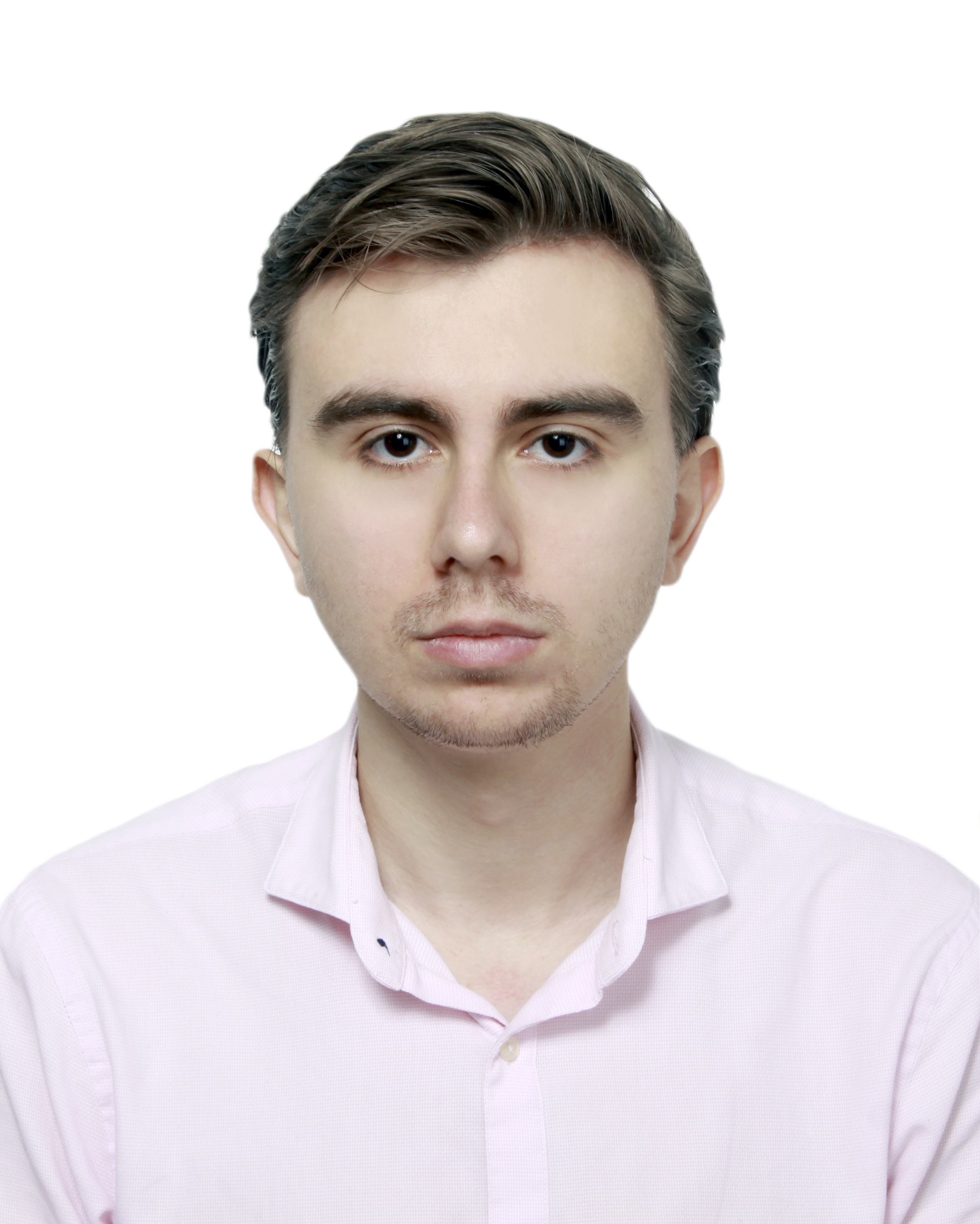}}]{Grigory Yaremenko}
	received his Bachelor’s degree from the faculty of Computational Mathematics and Cybernetic Faculty of Moscow State University in 2020, where the main focus of his studies and research were dynamical systems related disciplines and machine learning.
	His bachelor thesis was about an application of extended Kalman filtering to fluid dynamics and was awarded the first place at the Faculty's annual bachelor thesis competition.
	During his studies at Moscow State University, Grigory also worked at Vira Realtime as an industrial scientist and a software engineer from 2018 through 2020, producing data driven solutions for leak detection in oil pipelines.
	In 2022, Grigory was awarded a master's degree for completing Skoltech's "Data Science" program. During his time at Skoltech Grigory conducted theoretical research in stochastic stabilization and reinforcement learning.
\end{IEEEbiography}

\begin{IEEEbiography}
	[{\includegraphics[width=1in,height=1.25in,clip,keepaspectratio]{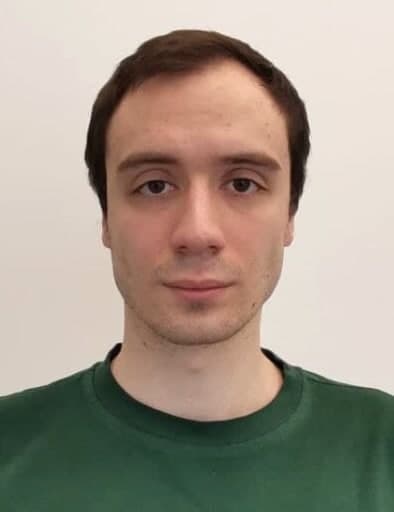}}]{Georgiy Malaniya}
	recieved a master's degree from the faculty of Mechanics and Mathematics of Moscow State University in 2021.
	His research has been mainly focused on control theory, symplectic geometry, algebraic topology and machine learning.
	Georgiy has industrial experience in research and development at Huawei Tech as a machine learning research engineer from 2019 through 2021 while developing time series models based on Bayesian methods.
	In 2021, Georgiy was granted a PhD position at the Skoltech's Artificial Intelligence in Dynamic Action Lab.
	Georgiy's areas of research include dynamical systems, reinforcement learning, manifold learning.
\end{IEEEbiography}
\vfill\null

%

\end{document}